\documentclass[11pt]{article}
\usepackage{etex}

\usepackage[margin=1in]{geometry}
\usepackage{stackengine}

\usepackage{stmaryrd}

\usepackage{tabularx}

\usepackage{color}
\usepackage[latin1]{inputenc}
\usepackage[T1]{fontenc}
\usepackage[normalem]{ulem}
\usepackage[english]{babel}
\usepackage{verbatim}
\usepackage{graphicx}
\usepackage{enumerate}
\usepackage{amsmath,amssymb,amsfonts,amsthm,amscd,mathrsfs}
\usepackage{array}
\usepackage{bbm}

\usepackage{dsfont}

\usepackage[T1]{fontenc}
\usepackage{babel}

\usepackage{url}

\usepackage{caption}
\usepackage{subcaption}
\usepackage[bookmarksopen, bookmarksnumbered]{hyperref}

\usepackage{cleveref}

\usepackage{tikz}
\usetikzlibrary{arrows}
\usetikzlibrary{arrows.meta}
\definecolor{wwhhii}{rgb}{1.,1.,1.}
\definecolor{rreedd}{rgb}{1.,0.,0.}
\definecolor{uuuuuu}{rgb}{0.26666666666666666,0.26666666666666666,0.26666666666666666}
\definecolor{darkgreen}{HTML}{0d8513}

\newtheorem{theorem}{Theorem}[section]

\newtheorem{lemma}[theorem]{Lemma}

\newtheorem{prop}[theorem]{Proposition}
\newtheorem{cor}[theorem]{Corollary}
\newtheorem{conj}[theorem]{Conjecture}

\theoremstyle{definition}

\newtheorem{rem}[theorem]{Remark}

\input xy
\xyoption{all}

\DeclareMathOperator{\erf}{erf}

\newcommand{\E}{\mathbb E}
\newcommand{\PP}{\mathbb P}

\newcommand{\Z}{\mathbb Z}
\newcommand{\R}{\mathbb R}

\newcommand{\N}{\mathbb N}

\newcommand{\cG}{\mathcal G}
\newcommand{\cF}{\mathcal F}
\newcommand{\cE}{\mathcal E}

\newcommand{\cH}{\mathcal H}

\makeatletter
\renewcommand\tableofcontents{%
  \null\hfill\textbf{\Large\contentsname}\hfill\null\par
  \@mkboth{\MakeUppercase\contentsname}{\MakeUppercase\contentsname}%
  \@starttoc{toc}%
}

\g@addto@macro\normalsize{%
  \setlength\abovedisplayskip{5pt}
  \setlength\belowdisplayskip{5pt}
  \setlength\abovedisplayshortskip{3pt}
  \setlength\belowdisplayshortskip{3pt}
}

\makeatother

\numberwithin{equation}{section}

\begin{document}

\title{Non-Markovianity of $2K-B$ and a degeneration}

\author{
Yang Chu
\thanks{Department of Statistics, University of California, Berkeley. e-mail: yang.chu@berkeley.edu}
\and
Lingfu Zhang
\thanks{Department of Statistics, University of California, Berkeley, and the Division of Physics, Mathematics and Astronomy, California Institute of Technology. e-mail: lingfuz@caltech.edu}
}
\date{}

\maketitle

\begin{abstract}
     We study the process $2K-B$, where $B$ is a standard 1-dimensional Brownian motion and $K$ is its concave majorant.
In light of Pitman's $2M-B$ theorem, it was recently conjectured by Ouaki and Pitman \cite{OP} that $2K-B$ has the law of the $\operatorname{BES}(5)$ process. The two processes share properties such as Brownian scaling, time inversion, and quadratic variation, as well as the same one-point distribution and infinitesimal generator. Yet it remains to be proved that $2K-B$ is Markovian.
However, we show that this conjecture is false.
To better understand the similarity between these two processes, we study a degeneration of $2K-B$. We show that it is a mixture of $\operatorname{BES}(3)$ and obtain other properties, including multiple-point distributions, the infinitesimal generator, and path decomposition at the future infimum.
We also further investigate the Markovian structure and the filtrations of $2K-B$, $B$, and $K$.
\end{abstract}

\begin{figure}[!ht] 
\centering
    \resizebox{0.6\hsize}{!}{
\includegraphics[]{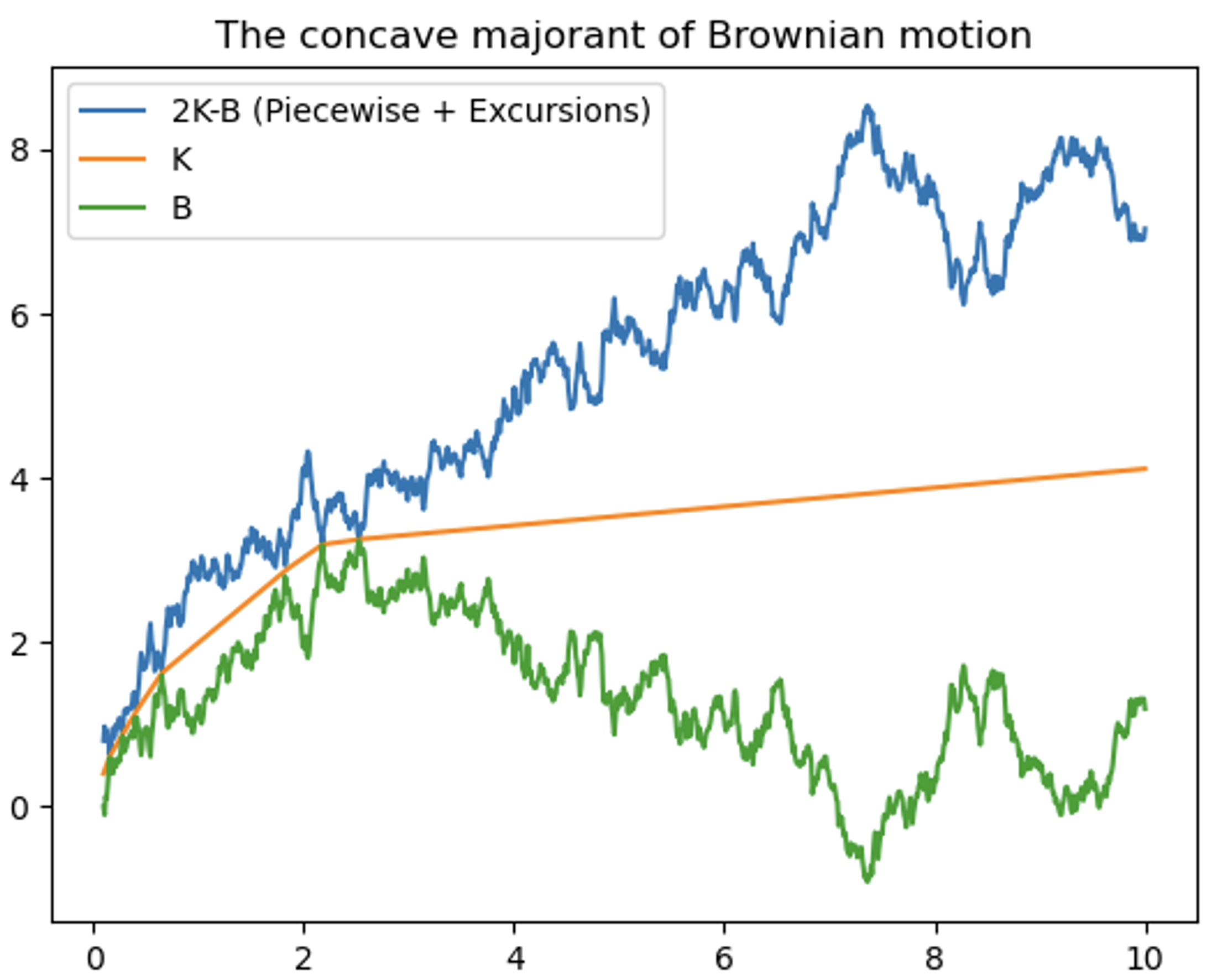}
}
\caption{A simulation of a Brownian motion $B$ and its concave majorant $K$, and the reflection $2K-B$.}\label{fig:kb}
\end{figure}

\section{Introduction} \label{sec:intro}

As one of the most fundamental random processes, Brownian motion (or Wiener process) has been intensively studied over the past century, with most of its properties now well understood.
Among them, a classical example is Pitman's $2M-B$ theorem: for a 1D Brownian motion $B:\R_{\ge 0}\to\R$ and its ``maximum process'' $M$, defined by $M(x):=\max_{0\le y\le x}B(y)$, the reflected process $2M-B$ has the same law as the Bessel process of dimension $3$ ($\operatorname{BES}(3)$), i.e., the radial part of a 3D Brownian motion \cite{Pit}.
Therefore, in some literature the function $x\mapsto \left(2\max_{0\le y \le x}f(y)\right) - f(x)$ (or $f(x)-\left(2\min_{0\le y \le x}f(y)\right)$) for a function $f:\R_{\ge 0}\to \R$ is known as the ``Pitman transform'' (see, e.g., \cite[Section 3.1]{DNV}).
Besides being elegant, the $2M-B$ theorem has stimulated much profound research, e.g., in Markov functions \cite{RogersPitman}, enlargement of filtration \cite{Jeulin_enlargement}, exponential functionals of Brownian motion \cite{matsumoto2000analogue, matsumoto2001analogue}, and random matrix theory \cite{o2002representation}.
Recently, the Pitman transform has played a crucial role in the study of models in the Kardar--Parisi--Zhang (KPZ) universality class.  
More specifically, a central family of models in this class are the \emph{last-passage percolation} models, such as the \emph{last passage} across a family of independent Brownian motions (see, e.g., \cite{WFS} for background).  
In \cite{DOV}, this model is shown to converge under KPZ scaling, and a key observation is that the \emph{last-passage value} is invariant under a ``melon'' operation, which generalizes the Pitman transform to multiple functions (see \cite[Section~3]{DNV}).  
This approach also relies on the fact that the ``melon'' operation maps independent Brownian motions to independent Brownian motions conditioned never to collide (see \cite[Theorem~2]{OYar}), which can be viewed as a generalization of the $2M-B$ theorem to multiple Brownian motions.

Besides reflecting against $M$, another natural operation for a random process is to take the concave majorant (resp.~convex minorant), which, for a function $f:U\to\R$ where $U\subset \R$, refers to the smallest (resp.~largest) function $c:U\to \R$ that is concave (resp.~convex) with $c\ge f$ (resp.~$c\le f$).
For a 1D Brownian motion $B$, it was shown by Groeneboom in \cite{Gro} to have an almost surely finite concave majorant $K$ on $\R_{\ge 0}$ (see \Cref{fig:kb}).
Since then, concave majorants and convex minorants of Brownian motion and other stochastic processes have been analyzed in many works; see, e.g., \cite{Pitman1983} for connecting Groeneboom's description of $K$ with the Williams path decomposition \cite{WDPath}, \cite{cin92} for connections between $K$ and processes arising in $M/M/\infty$ queueing, \cite{suidan} for analysis of the length of the longest segments of $K$, \cite{ber00} for convex minorants of the Cauchy process on a unit interval, \cite{abrampitm,bravpitm} for concave majorants of L\'evy processes, \cite{pitmanross} for convex minorants of Brownian motion and Brownian meanders on finite intervals, \cite{gonzalez2024convex} for L\'evy processes with smooth concave majorants and convex minorants, \cite{10.1214/20-EJP503} for convex minorants of L\'evy meanders, and \cite{OP} for a Markovian structure in $K$.  
There are also many applications, including isotonic regression in statistics \cite{isoto}, and ruin theory and mathematical finance \cite{gonzalez2022simulation, gonzalez2022geometrically}.  
See also \cite{abramevans,evansouaki} for studies of the related notion of Lipschitz minorants.

In the setting of a 1D Brownian motion $B$ and its concave majorant $K$ on $[0,\infty)$, the following conjecture was made by Ouaki and Pitman, by analogy with the $2M-B$ theorem.
\begin{conj}[\protect{\cite[Conjecture 1.3]{OP}}]  \label{coj:main}
The process $2K-B$ has the law of the $\operatorname{BES}(5)$ process.
\end{conj}
Here $\operatorname{BES}(5)$ denotes the Bessel process of dimension $5$, which is the radial part of a 5D Brownian motion.

\begin{figure}[t]
    \centering
\begin{subfigure}{0.48\hsize}
    \resizebox{0.99\hsize}{!}{
\includegraphics[]{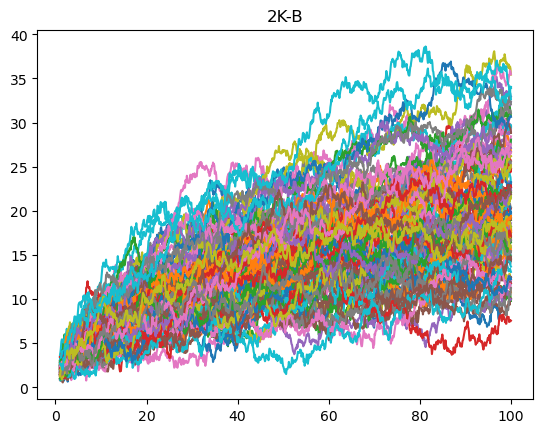}
    
}
\end{subfigure}
\begin{subfigure}{0.48\hsize}
    \resizebox{0.99\hsize}{!}{
\includegraphics[]{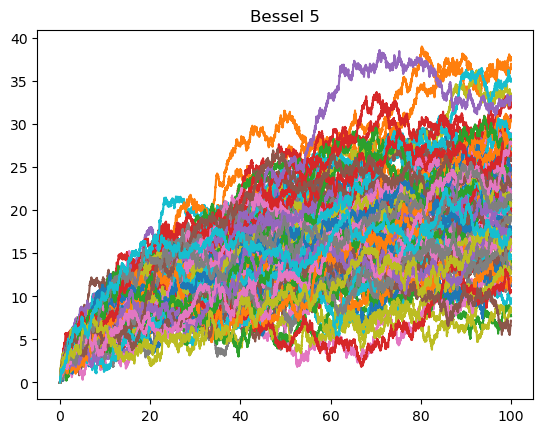}
}
\end{subfigure}
    \caption{100 samples of $2K-B$ and BES(5).}
    \label{fig:sim2}
\end{figure}

Much evidence supporting this conjecture is given in \cite{OP}: the processes $2K-B$ and $\operatorname{BES}(5)$ have the same one-point distribution and the same quadratic variation, and both are invariant under Brownian scaling and time inversion. 
Numerical simulations also seem to favor this conjecture (see \Cref{fig:sim2} and \Cref{sec:simu}).
Moreover, it is even shown that $2K-B$ has the same ``infinitesimal generator'' as $\operatorname{BES}(5)$.
\begin{prop}[\cite{OP}, Proposition 4.1] \label{OP_Prop_infi}
Let $\varphi$ be a twice continuously differentiable function with compact support on $(0, \infty)$, and let $t, z>0$. Then
\begin{multline*}
\lim_{s\to 0_+} s^{-1} \E[\varphi(2K(t+s)-B(t+s)) - \varphi(2K(t)-B(t)) \mid 2K(t)-B(t)=z]\\=
2z^{-1}\varphi'(z) + \varphi''(z)/2.
\end{multline*}
\end{prop}
\begin{rem} \label{OP_Rem_infi}
We note that \cite[Proposition 4.1]{OP} is stated only for $t=1$, and the general case $t>0$ follows from the invariance of $B$ and $K$ under diffusive rescaling.
In particular, we see that this ``infinitesimal generator'' does not depend on $t$. If $2K-B$ were Markovian, it would be a time-homogeneous Markov process.

For $R$ being the $\operatorname{BES}(5)$ process, it is time-homogeneous Markovian, with 
\[
\lim_{s\to 0_+} s^{-1} \E[\varphi(R(t+s)) - \varphi(R(t)) \mid R(t)=z] =
2z^{-1}\varphi'(z) + \varphi''(z)/2.
\]
See, e.g., \cite[Chapter XI, Section 1]{revuz_yor}.
Therefore, if $2K-B$ were shown to be Markovian, it would have to be $\operatorname{BES}(5)$, and Conjecture~\ref{coj:main} would follow.
\end{rem}

A Markovian structure of $K$ and $B$ is given in \cite{OP}.
More specifically, it is known (since \cite{Gro,Pitman1983}, quoted as \Cref{Pitman_seq} below) that $K$ is increasing and piecewise linear, with infinitely many linear segments that accumulate only at zero and at infinity; and knowing $K$, the process $K-B$ is a concatenation of independent Brownian excursions between \emph{vertices}, where $t$ is a vertex if $K(t)=B(t)$.
It is then shown in \cite{OP} that $\Psi = (K', K, K-B, H)$ is a time-homogeneous Markov process in $\R^4$, where $K'$ is the right-hand derivative of $K$, and $H(t)= D(t) - t$, with $D(t)$ denoting the next vertex time after $t$ (see \Cref{OP_THM_markov} below).

As pointed out at the beginning of \cite[Section 4]{OP}, based on the above Markovian analysis, the hope was to use the theory of Markov functions to prove that $2K-B$ is Markovian by checking the intertwining condition (see \cite{RogersPitman}), thereby proving Conjecture~\ref{coj:main}.
However, we show that this is not possible.
\begin{theorem}  \label{thm:nb5}
The process $2K-B$ is not Markovian, and therefore not equal in distribution to the $\operatorname{BES}(5)$ process.
\end{theorem}
A short proof by contradiction is given in \Cref{sec:short}.

As this disproves Conjecture~\ref{coj:main}, we further investigate the process $2K-B$ from several perspectives. 

In \Cref{sec:limit}, we consider a degeneration of $2K-B$, which can be interpreted as ``conditioning on flat $K'$''. More precisely, for any $s, z > 0$, let $Z_{s,z}$ be the process with the law of $\bigl(2K(s+t) - B(s+t)\bigr)_{t \ge 0}$ knowing $2K(s) - B(s) = z$. For fixed $z > 0$, the process $Z_{s,z}$ converges weakly to a process $Z_z$ as $s \to \infty$. This convergence will be rigorously proved in Proposition~\ref{limZdef}, where we show that $Z_{s,z}$ is well-defined and prove weak convergence as $s \to \infty$ in the compact-open topology.

As $z > 0$ is typically fixed, we shall usually omit it from the subscript and write $Z = Z_z$. There are several motivations for studying this degeneration $Z$. First, the analysis of $Z$ provides several alternative proofs of \Cref{thm:nb5}. For instance, if $2K-B$ were $\operatorname{BES}(5)$, the law of $Z_{s,z}$ would depend only on $z$, not on the time parameter $s$. Thus, one would expect the limit $Z$ to be $\operatorname{BES}(5)$ starting at $z$ if Conjecture~\ref{coj:main} holds. Showing that $Z$ is not $\operatorname{BES}(5)$ therefore yields another proof of \Cref{thm:nb5}. Second, $Z$ inherits many properties from $2K-B$, and its study helps explain why $2K-B$ resembles $\operatorname{BES}(5)$ so closely. Finally, $Z$ is more tractable than $2K-B$, and analyzing it sheds light on aspects of $2K-B$ that are otherwise difficult to access.

Our results on $Z$ include 
\begin{enumerate}
    \item showing that $Z$ is essentially a mixture of $\operatorname{BES}(3)$ (Proposition~\ref{BESmix}), by which we mean a $\operatorname{BES}(3)$ with a random starting point, together with a path decomposition at the future infimum (Proposition~\ref{inf_path_decom});
    \item computing its multi-point distribution (Proposition~\ref{Ztfinev});
    \item analyzing its infinitesimal generator in the context of a general mixture of $\operatorname{BES}(3)$ (Proposition~\ref{p:infgen});
    \item establishing a future infimum reflection property (Proposition~\ref{futinf}).
\end{enumerate}

In \Cref{sec:fil}, we further analyze the Markovian structure discovered in \cite{OP}.
We systematically study other natural projections of $(K', K, K-B, H)$ and determine whether each of them is Markovian.
One conclusion is that none of the natural projections of $(K', K, K-B, H)$ that contain $2K-B$ is Markovian.

In \Cref{sec:fX}, we analyze the filtration of $2K-B$, in comparison with the filtrations of $K$ and $B$.
In particular, we establish a semimartingale decomposition of $2K-B$ in an enlarged filtration.

\subsection*{Setup and notation}
We define some notations that are widely used throughout.

\medskip

\noindent\textbf{Basic notation and distribution.} For a random variable $U$ taking values in $\mathbb{R}^d$ with a probability density with respect to Lebesgue measure on $\mathbb{R}^d$, we denote that probability density function by $f_U$. Similarly, for two random variables $U$ and $V$ both taking values in Euclidean spaces (not necessarily of the same dimension), the notation $f_{U \mid V=v}(\cdot)$ is used to denote a regular conditional probability density of $U$ given $V=v$.

For $a, b>0$, we let $\beta_{a, b}$ denote a random variable with the $\operatorname{Beta}(a, b)$ density
$$
f_{\beta_{a, b}}(u)=\frac{\Gamma(a+b)}{\Gamma(a) \Gamma(b)} u^{a-1}(1-u)^{b-1} \mathds{1}_{\{0<u<1\}},
$$
where $\Gamma$ is the Gamma function. 

We let $\phi$ denote the probability density of a standard normal random variable:
$$
\phi(x)=\frac{1}{\sqrt{2 \pi}} \exp \left(-\frac{x^2}{2}\right).
$$

For any random variable $U$, we denote by $\sigma(U)$ the sigma-algebra generated by it.

We use the notation $\wedge$ to denote the minimum of two real numbers, and $\vee$ to denote the maximum of two real numbers, or the join of two sigma-algebras.

\medskip

\noindent\textbf{Random processes.}
We always use $B$ to denote a Brownian motion with unit diffusivity. For $\mu \in \mathbb{R}$, we use the notation $B_\mu(t):=B(t)+\mu t$ for a Brownian motion with drift $\mu$. For each level $h>0$ and drift $\mu \ge 0$, define the first and last passage times
$$
T_{\mu, h}:=\inf \left\{t>0: B_\mu(t)=h\right\} \quad \text{and} \quad G_{\mu, h}:=\sup \left\{t>0: B_\mu(t)=h\right\}.
$$
Then $T_{\mu, h}$ has the inverse Gaussian distribution with parameters $(\mu, h)$, whereas $G_{\mu, h}$ has the size-biased inverse Gaussian distribution with parameters $(\mu, h)$. Their respective densities are denoted by $f_{\mu, h}$ and $f_{\mu, h}^*$, and are given by (see, e.g., \cite[Notation 1.1]{OP})
\begin{equation}\label{inverse_gaussian}
f_{\mu, h}(t):=f_{T_{\mu, h}}(t)=\frac{h}{\sqrt{2 \pi t^3}} \exp \left(-\frac{(h-\mu t)^2}{2 t}\right) \mathds{1}_{t>0}, \quad
f_{\mu, h}^*(t):=f_{G_{\mu, h}}(t)=\frac{\mu}{h} t f_{\mu, h}(t).
\end{equation}
Note that when $\mu \to 0$, $T_{\mu, h} \to T_{0, h}$, which is the $\frac{1}{2}$-stable L\'evy distribution with density
\[
f_{0, h}(t) = \frac{h}{\sqrt{2 \pi t^3}} \exp \left(-\frac{h^2}{2 t}\right) \mathds{1}_{t>0},
\]
whereas $G_{\mu, h} \to \infty$ almost surely, since Brownian motion without drift does not have a finite last passage time.

For any $h, \mu \ge 0$ and $n \in \N$, the $\operatorname{BES}_h(n, \mu)$ process is the \emph{$n$-dimensional Bessel process} starting at $h$ with drift $\mu$.
Namely, it is the radial part of $n$-dimensional Brownian motion with a drift vector whose norm is $\mu$. When either $h$ or $\mu$ equals $0$, we omit them from the notation.

By an \emph{$n$-dimensional Bessel bridge} from $(a,h)$ to $(b,g)$, for some $a<b$ and $g\ge 0$, we refer to $t\mapsto W(t-a)$, where $W$ is $\operatorname{BES}_h(n)$ restricted to $[0,b-a]$, conditional on taking the value $g$ at time $b-a$ (see, e.g., \cite[Chapter XI, Section 3]{revuz_yor}).
By a \emph{Brownian excursion} on the time interval $[a, b]$, we refer to a 3-dimensional Bessel bridge from $(a,0)$ to $(b,0)$ (see \cite{williams1970decomposing}).
We refer interested readers to \cite[Chapters~3,~6,~and~11]{revuz_yor} for further information on these processes.

For any random process $W$ on $[0,\infty)$, we denote by $(\cF^W_t)_{t\ge 0}$ the filtration generated by $W$.

\medskip

\noindent\textbf{Concave majorant.} 
We use $K$ to denote the concave majorant of a Brownian motion $B$ on $[0,\infty)$. We also write $X = 2K - B$.

Recall that \emph{vertices} are those $t>0$ such that $K(t)=B(t)$.
For each $t>0$, let $G(t)$ and $D(t)$ be the left-hand and right-hand vertices of the segment of $K$ straddling $t$. Then we almost surely have $0 < G(t) < t < D(t)$. We also define $I(t) = K(t) - t K^{\prime}(t)$ to be the intercept at $0$ of the line extending this segment. 
For any $a>0$, let
\[
\theta(a) = \sup \left\{t>0: B(t)-a t = \sup_{u>0}(B(u)-a u)\right\} = \inf\{t>0: K'(t)<a\}.
\]
This is the first time at which the slope of the concave majorant drops below $a$.
Then the set of vertices is precisely $\{\theta(a): a>0\}$, the set of jump times of $K'$.

\section{Preliminary on the concave majorant of Brownian motion}
As a starting point of our discussion, we state a sequential description of the joint law of $K$ and $B$ from \cite{Pitman1983,Gro}, in terms of the vertices and segment lengths.

Fix a number $b > 0$, and we enumerate the vertices starting from $\theta(b)$. Namely, let $V_0 := \theta(b)$, and let $V_i$ (resp.~$V_{-i}$) be the $i$-th vertex to the right (resp.~left) of $V_0$ for any $i \in \N$.
For any $i \in \Z$, denote
\begin{align*}
& T_i = V_{i+1} - V_i, \\
& \alpha_i = \frac{K(V_{i+1}) - K(V_i)}{T_i}.
\end{align*}
Note that we then have $V_i = \sum_{j<i} T_j$ and $K(V_i) = \sum_{j<i} \alpha_j T_j$. Therefore, the process $K$ is determined by the sequence $\{T_i, \alpha_i\}_{i \in \mathbb{Z}}$.

\begin{theorem}[\cite{Gro,Pitman1983}] \label{Pitman_seq} 
The law of $\{T_i, \alpha_i\}_{i \in \mathbb{Z}}$ can be described as follows.
\begin{enumerate}
    \item $\alpha_0 \sim \operatorname{Unif}(0,b)$; and for each $i \in \N$, conditional on $\alpha_0, \ldots, \alpha_{i-1}$, the law of $\alpha_i$ is $\operatorname{Unif}(0, \alpha_{i-1})$;
    \item $\alpha_{-1}$ has probability density $b x^{-2}$ on $(b, \infty)$; and for each $i \in \N$, $i \ge 2$, conditional on $\alpha_{-1}, \ldots, \alpha_{-i+1}$, the law of $\alpha_{-i}$ has probability density $\alpha_{-i+1} x^{-2}$ on the interval $(\alpha_{-i+1}, \infty)$;
    \item the sequences $\{\alpha_i\}_{i=0}^\infty$ and $\{\alpha_{-i}\}_{i=1}^\infty$ are independent;
    \item conditional on all the slopes $\{\alpha_i\}_{i \in \Z}$, the lengths $T_i$ of the segments are independent and have the $\operatorname{Gamma}\!\left(\frac{1}{2}, \frac{1}{2} \alpha_i^2\right)$ distribution.
\end{enumerate}
Moreover, given $\{T_i, \alpha_i\}_{i \in \mathbb{Z}}$, the processes $K-B$ on each $[T_{i-1}, T_i]$ are independent Brownian excursions.
\end{theorem}

We now recall two results from \cite{OP} on a Markovian structure involving $K$.
\begin{prop}[\cite{OP}, Proposition 1.1] \label{OP_Prop_density}
The density of $\left(K^{\prime}(1), I(1), K(1)-B(1), \frac{1}{G(1)}, D(1)\right)$ is\footnote{We note that in \cite{OP}, there is a typo: the factor of $y^2$ in the expression for $f_5(a, b, y, v, w)$ is omitted.}
$$
\begin{aligned}
f_5(a, b, y, v, w)
&= \sqrt{\frac{2}{\pi^3(v-1)^3(w-1)^3}}\, a b(w v - 1) y^2 \\
&\quad \times \exp\!\left(-\frac{1}{2}\left(b^2 w + 2 a b + a^2 v + y^2 \frac{w v - 1}{(v-1)(w-1)}\right)\right)
\mathds{1}_{\{w, v>1\}} \mathds{1}_{\{a, b, y>0\}}.
\end{aligned}
$$
In particular, the following marginals take simpler forms:
\begin{enumerate}
    \item The joint density of $\left(K^{\prime}(1), I(1), K(1)-B(1)\right)$ at $(a, b, y)$ is
    $$
    f_3(a, b, y)=4 y(a+b+y)\,\phi(a+b+y)\,\mathds{1}_{\{a, b, y>0\}}.
    $$
    \item The conditional density of $D(1)-1$ given $K^{\prime}(1)=a$, $I(1)=b$, and $K(1)-B(1)=y$ is
    (recall the inverse Gaussian and size-biased inverse Gaussian distributions introduced in \eqref{inverse_gaussian})
    $$
    \frac{a}{a+b+y}\, f_{a, y} + \frac{b+y}{a+b+y}\, f_{a, y}^*.
    $$
\end{enumerate}
\end{prop}

\begin{theorem}[\cite{OP}, Theorem 1.2] \label{OP_THM_markov}
The process
$$
\Psi(t):=\left(K^{\prime}(t), K(t), K(t)-B(t), D(t)-t\right), \quad t \geq 0,
$$
is a time-homogeneous Markov process. Its semigroup $P_t((a, k, y, w), \cdot)$ is given as follows:
\begin{enumerate}
    \item If $t<w$, then $P_t((a, k, y, w), \cdot)$ has the law of
    $$
    (a,\, k+a t,\, R(t),\, w-t),
    $$
    where $R$ is a 3-dimensional Bessel bridge from $(0, y)$ to $(w, 0)$.
    \item If $t \geq w$, then $P_t((a, k, y, w), \cdot)$ has the law of
    $$
    \left(a-C^{\prime}(t-w),\, k+a t-C(t-w),\, R(t-w)-C(t-w),\, D^R(t-w)-(t-w)\right),
    $$
    where $R$ is a $\operatorname{BES}(3, a)$ process, $C$ is its convex minorant, and $D^R(t-w)$ is the first vertex of $C$ after time $t-w$.
\end{enumerate}
\end{theorem}

\section{A proof of \Cref{thm:nb5} by contradiction}   \label{sec:short}
We now give a straightforward disproof of Conjecture~\ref{coj:main}.
We start with the following preparatory lemmas. 

\begin{lemma}  \label{lem:sdebb}
Take $h \ge 0$ and an integer $n \ge 2$, and let $W:[0,1] \to \R_{\ge 0}$ be an $n$-dimensional Bessel bridge from $(0,h)$ to $(1,0)$.
Then $W$ satisfies the following stochastic differential equation:
\[
W(t) - h = B(t) + \int_0^t \left( \frac{n-1}{2W(s)} - \frac{W(s)}{1-s} \right) ds,
\]
for any $t \in [0, 1]$, where $B:[0,\infty) \to \R$ is a Brownian motion with unit diffusivity.
\end{lemma}

In the case where $h=0$, this statement is given in \cite[Chapter XI, Exercise 3.11]{revuz_yor}, with a sketched proof that also applies to general $h>0$.
For the reader's convenience, we reproduce the proof here.
\begin{proof}[Proof of Lemma \ref{lem:sdebb}]
For any $t, x, y > 0$, let $p_n(t, x, y)$ denote the transition density of the $n$-dimensional Bessel process from $x$ to $y$ in time $t$.
Then, for any $t \in (0, 1)$, the Radon--Nikodym derivative of $W$ on $[0, t]$ with respect to the $\operatorname{BES}_h(n)$ process on $[0, t]$ is
\[
 \lim_{y \to 0}\frac{p_n(1-t, W(t), y)}{p_n(1, h, y)} 
 = (1-t)^{-n/2} \exp\left(\frac{h^2}{2} - \frac{W(t)^2}{2(1-t)}\right).
\]
This equality follows from the fact that, for any $t, x > 0$, the limit $y \to 0$ of $p_n(t, x, y)\,y^{-(n-1)}$ equals (up to a universal constant) the transition density of an $n$-dimensional Brownian motion from a point on the sphere of radius $x$ to the origin in time $t$.

On the other hand, let $\tilde{W}$ denote the $\operatorname{BES}_h(n)$ process. It satisfies the stochastic differential equation (see, e.g., \cite[Chapter XI, (1.26)]{revuz_yor})
\[
\tilde{W}(t) - h = \tilde{B}(t) + \int_0^t \frac{n-1}{2\tilde{W}(s)}\, ds,
\]
for any $t \in [0, 1]$, where $\tilde{B}:[0,\infty) \to \R$ is a Brownian motion with unit diffusivity.
The result then follows from Girsanov's theorem (as stated in \cite[Chapter VIII, Theorem 1.4]{revuz_yor}), applied to the laws of $W$ and $\tilde{W}$ on $[0, t]$ for any $t \in (0,1)$.
\end{proof}



\begin{lemma}  \label{lem:pre-coup}
Take any interval $[a, b]$, with $0 \le h_1 \le h_2$, and let $n \ge 2$ be an integer.
Let $W_1:[a,b] \to \R_{\ge 0}$ be an $n$-dimensional Bessel bridge from $(a,h_1)$ to $(b,0)$, and let $W_2:[a,\infty) \to \R_{\ge 0}$ be an $n$-dimensional Bessel process starting at $W_2(a)=h_2$. Then we can couple them on the same probability space so that, almost surely, $W_1(t)\le W_2(t)$ for all $t \in [a, b]$. 
\end{lemma}

\begin{proof}
By Lemma~\ref{lem:sdebb}, the process $W_1$ satisfies the stochastic differential equation
\[
W_1(t) - h_1 = B(t-a) + \int_a^t \left( \frac{n-1}{2W_1(s)} - \frac{W_1(s)}{b-s} \right) ds,
\]
for any $t \in [a, b]$.
On the other hand, $W_2$ satisfies the stochastic differential equation (see, e.g., \cite[Chapter XI, (1.26)]{revuz_yor})
\[
W_2(t) - h_2 = B(t-a) + \int_a^t \frac{n-1}{2W_2(s)}\, ds,
\]
for any $t \in [a, b]$.
By coupling $W_1$ and $W_2$ with the same driving Brownian motion $B$, we obtain
\[
(W_2(t)-h_2)-(W_1(t)-h_1) = \int_a^t \left( (W_1(s)-W_2(s))\frac{n-1}{2W_1(s)W_2(s)} + \frac{W_1(s)}{b-s} \right) ds.
\]
Suppose that the set $\{t \in [a, b] : W_1(t) > W_2(t)\}$ is nonempty. Since $a$ is not in this set, there exist $a \le t_- < t_+ \le b$ such that $W_1(t_-) = W_2(t_-)$, while $W_1(t) > W_2(t)$ for all $t \in (t_-, t_+]$.
However, we have
\[
W_2(t_+) - W_1(t_+) = \int_{t_-}^{t_+} \left( (W_1(s)-W_2(s))\frac{n-1}{2W_1(s)W_2(s)} + \frac{W_1(s)}{b-s} \right) ds.
\]
The left-hand side is $<0$, while the right-hand side is $>0$, yielding a contradiction.
\end{proof}

For simplicity of notation, in the rest of this section we use $R_{n,h}$ to denote a $\operatorname{BES}_h(n)$ process.
The main task is to show that $X$ and $R_{5,0}$ are not equal in distribution.

\begin{lemma}   \label{lem:Xdom}
For any $t>0$, we can couple $X-K$ on $[t, \infty)$, conditional on $X(t)-K(t)$, with $R_{3,X(t)-K(t)}$, such that almost surely $X(t+s)-K(t+s)\le R_{3,X(t)-K(t)}(s)$ for all $s\ge 0$.
\end{lemma}

\begin{proof}
We let $t_0 = t$, and let $t_1 < t_2 < \cdots$ be the sequence of all vertices in $(t, \infty)$.
By \Cref{Pitman_seq}, conditional on the sequence $t_1 < t_2 < \cdots$, the law of $X - K = K - B$ on $[t_0, \infty)$ is as follows.
On $[t_0, t_1]$, it is a 3-dimensional Bessel bridge from $(t_0, X(t_0) - K(t_0))$ to $(t_1, X(t_1) - K(t_1)) = (t_1, 0)$. 
For each $i \in \N$, $X - K$ is a Brownian excursion on $[t_i, t_{i+1}]$
(recall that a Brownian excursion on the time interval $[a, b]$ is a 3-dimensional Bessel bridge from $(a, 0)$ to $(b, 0)$, as shown in \cite{williams1970decomposing}).
All these Bessel bridges and Brownian excursions are independent. 

We then couple $R_{3, X(t) - K(t)}$ with $X - K$ inductively, conditional on the sequence $t_1 < t_2 < \cdots$.
First, by Lemma~\ref{lem:pre-coup}, we couple $R_{3, X(t) - K(t)}$ on $[0, t_1 - t_0]$ and $X - K$ on $[t_0, t_1]$ so that 
\[
X(t + s) - K(t + s) \le R_{3, X(t) - K(t)}(s), \quad \text{for all } s \in [0, t_1 - t_0].
\]
Next, for any $i \in \N$, assume that we have a coupling between $R_{3, X(t) - K(t)}$ on $[0, t_i - t_0]$ and $X - K$ on $[t_0, t_i]$, such that 
\[
X(t + s) - K(t + s) \le R_{3, X(t) - K(t)}(s), \quad \text{for all } s \in [0, t_i - t_0].
\]
Note that for $R_{3, X(t) - K(t)}$, since it is Markovian, conditional on its path over $[0, t_i - t_0]$, it is a 3-dimensional Bessel process starting from $R_{3, X(t) - K(t)}(t_i - t_0)$ on $[t_i - t_0, \infty)$.
Therefore, we can apply Lemma~\ref{lem:pre-coup} again to couple $R_{3, X(t) - K(t)}$ on $[t_i - t_0, t_{i+1} - t_0]$ and $X - K$ on $[t_i, t_{i+1}]$ so that 
\[
X(t + s) - K(t + s) \le R_{3, X(t) - K(t)}(s), \quad \text{for all } s \in [t_i - t_0, t_{i+1} - t_0].
\]
Thus, the conclusion follows.
\end{proof}

\begin{proof}[Proof of \Cref{thm:nb5}]
We argue by contradiction and assume that $X$ has the same law as $R_{5,0}$.

Take a small enough $\epsilon > 0$ (e.g., $\epsilon < 0.1$ suffices), and let $\cE$ be the following event:
\[
X(1) < \epsilon, \quad X(2) \in [1, 1 + \epsilon].
\]
Let $\cE'$ denote the same event with $X$ replaced by $R_{5,0}$.
By the above assumption, we have $\PP[\cE] = \PP[\cE'] > 0$. 

The reason for considering these events is that conditioning on them highlights the difference between $X = 2K - B$ and $R_{5,0}$.
Namely, the condition $X(1) < \varepsilon$ in the event $\cE$ forces the slope of the majorant to be nearly flat, so the instantaneous drift of $X$ at time $2$ (assuming $X(2) \in [1, 1 + \epsilon]$) behaves like that of a $\operatorname{BES}(3)$ process, which is smaller than the drift of a $\operatorname{BES}(5)$ process.
By contrast, since $R_{5,0}$ is Markovian, its drift at time $2$, conditional on $R_{5,0}(2)$, is independent of $\cE'$ and remains that of a $\operatorname{BES}(5)$ process.

We now justify this rigorously.
We consider the drift of $X$ at time $2$, conditional on $\cE$.
Since $K \le X$, we must have $K(1) < \epsilon$ under $\cE$.
This implies that $0 < K(t) - K(s) < \epsilon (t-s)$ for any $1 \le s < t$, since $K$ is concave and increasing.
Then we have
\begin{align*}
\limsup_{s\to 0_+} s^{-1}\E[X(2+s)-X(2)&\mid \cE]
\le \limsup_{s\to 0_+} s^{-1}\E[K(2+s)-K(2)\mid \cE]\\
+ & \limsup_{s\to 0_+} s^{-1}\E[(X(2+s)-K(2+s))-(X(2)-K(2))\mid \cE]\\
&\le \epsilon + \lim_{s\to 0_+} s^{-1}\E[R_{3,X(2)-K(2)}(s)-(X(2)-K(2)) \mid \cE]\\
&= \epsilon + \lim_{s\to 0_+} s^{-1}\E[R_{3,X(2)-K(2)}(s)-R_{3,X(2)-K(2)}(0) \mid \cE].
\end{align*}
In the second inequality, we used the fact that $K(2+s)-K(2) < \epsilon s$ under $\cE$ (to bound the first term) and Lemma~\ref{lem:Xdom} (to bound the second term).
Then, by the stochastic differential equation satisfied by the 3-dimensional Bessel process $R_{3, X(2) - K(2)}$ (see, e.g., \cite[Chapter XI, (1.26)]{revuz_yor}), the limit in the last line equals $\E\left[\frac{1}{X(2) - K(2)} \mid \mathcal{E} \right]$, which is at most $\frac{1}{1 - 2\epsilon}$.
Therefore, we have
\[
\limsup_{s\to 0_+} s^{-1}\E[X(2+s)-X(2)\mid \cE] \le \epsilon + \frac{1}{1 - 2\epsilon}.
\]
On the other hand,
\[
\lim_{s\to 0_+} s^{-1}\E[R_{5,0}(2+s)-R_{5,0}(2)\mid \cE'] = \E\left[\frac{2}{R_{5,0}(2)} \,\middle|\, \cE'\right] \ge \frac{2}{1 + \epsilon}.
\]
By taking $\epsilon$ small (e.g., $\epsilon < 0.1$), we obtain a contradiction.

Since $X$ is not equal in law to $R_{5,0}$, it follows from Proposition~\ref{OP_Prop_infi} and Remark~\ref{OP_Rem_infi} that it is not Markovian.
\end{proof}

\section{A flat degeneration of $X$}\label{sec:limit}

In this section we fix $z>0$ and study the process $Z_z = Z$, which was described after \Cref{thm:nb5} and will be rigorously constructed below in Proposition~\ref{limZdef}.
This process can be viewed as a degeneration of $X$, as will be explained shortly. 
On the other hand, our initial motivation for considering $Z$ comes from Pitman's $2M-B$ theorem: there, a key property is that for any $t>0$, the law of $M(t)$ conditional on $2M(t)-B(t)=z$ is $\operatorname{Unif}(0,z)$. 
In particular, this law is independent of $t$, leading to the fact that $2M-B$ is time-homogeneous Markovian.
It is therefore natural to consider $\Psi(t)$ (from \Cref{OP_THM_markov}) conditional on $X(t)=z$.
Below we record this conditional distribution, which does depend on $t$, indicating that $X$ is not time-homogeneous Markovian.

\begin{lemma} \label{condi-law}
Conditional on $X(t)=z$, we have:
\begin{enumerate}
    \item The density of $K(t)-B(t)$ at $y$ is $\frac{6 y(z-y)}{z^3} \mathds{1}_{0<y<z}$ (thus it does not depend on $t$).\\
    Consequently, since $K(t)=X(t)-(K(t)-B(t))$, the density of $K(t)$ at $y$ is also $\frac{6 y(z-y)}{z^3} \mathds{1}_{0<y<z}$.
    \item When further conditioning on $K(t)$, the law of $K'(t)$ is $\operatorname{Unif}(0,K(t)/t)$.
    \item When further conditioning on $K'(t)=a$ and $K(t)=y$, the density of $D(t)-t$ is given by $\frac{a t}{z} f_{a, z-y} + \left(1 - \frac{a t}{z}\right) f_{a, z-y}^*$.
\end{enumerate}
\end{lemma}

\begin{proof}
By the scaling property of the process, 
\[
\left(K'(t), K(t), K(t)-B(t), D(t)-t\right)
\stackrel{d}{=}
\left(\frac{1}{\sqrt{t}} K'(1), \sqrt{t} K(1), \sqrt{t}(K(1)-B(1)), t(D(1)-1)\right),
\]
and the result follows immediately from Proposition~\ref{OP_Prop_density}.
\end{proof}
From here, we see that the dependence of the law on time $t$ arises solely through the slope $K'(t)$. (While $a t = K'(t)t$ appears in the density of $D(t)-t$, the random variable $t K'(t)$ is $\operatorname{Unif}(0, K(t))$, which is independent of $t$.) 
To ``remove this dependence,'' we fix $z$ and let time tend to infinity, obtaining the process $Z$. This forces $K' \to 0$. 

Below we study various aspects of $Z$, as announced after \Cref{thm:nb5}.
Along the way, we provide several alternative proofs of \Cref{thm:nb5}.
This is because, if $X$ were a $\operatorname{BES}(5)$ process (and hence Markovian), $Z$ would be $\operatorname{BES}_z(5)$; however, this can be refuted by various properties of $Z$, to be revealed in this section.
Our study of $Z$ also aims to better understand the original process $X$, and in particular, why it resembles $\operatorname{BES}(5)$ in many aspects.

\subsection{Description of the limiting process}

\begin{prop}  \label{limZdef}
For fixed $z>0$, the law of $(X(s+t)\mid X(s)=z)_{t\ge 0}$ converges (in the topology of uniform convergence on each compact
subset of $[0,\infty)$) as $s \to \infty$, and the limiting process $(Z(t))_{t\ge 0}$ can be described as follows:
\begin{enumerate}
    \item Let $\tilde{K}$ be a random variable with density $\frac{6 k(z-k)}{z^3} \mathds{1}_{\{0<k<z\}}$, and define $P = \frac{U\tilde{K}}{z}$, where $U$ is uniform on $[0, 1]$ and independent of $\tilde{K}$.
    \item Given $\tilde{K}$ and $P$, with probability $1-P$, let $Z = \tilde{K} + Y$, where $Y$ is $\operatorname{BES}_{z-\tilde{K}}(3)$; with probability $P$, let $Z = z + B$ on $[0, \tau]$, where $B$ is a Brownian motion starting at $0$, $\tau$ is the first time $B$ hits $\tilde{K} - z$, and $(Z(\tau+t)-Z(\tau))_{t\ge 0}$ is a $\operatorname{BES}(3)$ process independent of $B$.
\end{enumerate}
\end{prop}

\begin{proof}
From Lemma~\ref{condi-law}, parts~1 and~2, as $s \to \infty$, $K(s)$ and $sK'(s)$ conditional on $X(s)=z$ converge to $\tilde K$ and $zP$, respectively.
Then, by Lemma~\ref{condi-law}, part~3, the $s \to \infty$ limit of the law of $D(s)-s$ conditional on $X(s)=z$ and $K'(s), K(s)$ is as follows: it is a mixture of $\infty$ (with probability $1-P$) and $f_{0,z-\tilde K}$ (with probability $P$).

Note that by \Cref{OP_THM_markov}, conditional on $K'(s)$, $K(s)$, $B(s)$, and $D(s)$,
the process $K-B = X-K$ on $[s, D(s)]$ is a 3-dimensional Bessel bridge from $(s, X(s)-K(s))$ to $(D(s), 0)$, and on $[D(s), \infty)$ it is given by $\operatorname{BES}(3, K'(s))$.
Then, as $s \to \infty$, $(X(s+t))_{t \ge 0} - \tilde K$, conditional on $X(s)=z$, converges to the following process:  
with probability $1-P$, take $\operatorname{BES}_{z-\tilde K}(3)$; with probability $P$, take a
3-dimensional Bessel bridge from $(0, z-\tilde K)$ to $(\tau_*, 0)$ concatenated with a $\operatorname{BES}(3)$ process, where $\tau_*$ is a random variable with density $f_{0,z-\tilde K}$.

On the other hand, for $B$ and $\tau$ in the statement, the probability density of $\tau$ is also given by $f_{0,z-\tilde K}$, and conditional on $\tau$, $B+z-\tilde K$ on $[0, \tau]$ is precisely a 3-dimensional Bessel bridge from $(0, z-\tilde K)$ to $(\tau, 0)$.
Therefore, the conclusion follows.
\end{proof}

We next give two alternative descriptions of the limiting process $Z$.
For these, we need the probability density function of the infimum of Bessel processes, which is classical and can be derived from martingale arguments (see, e.g., \cite{hamana2013probability}).

\begin{lemma}  \label{infimum}
For any $h>0$ and integer $n \ge 3$, consider the minimum of $\operatorname{BES}_h(n)$: its probability density at $x$ is $(n-2)x^{n-3}h^{-(n-2)}\mathds{1}_{0<x<h}$.
\end{lemma}

Another result we need is the Williams path decomposition \cite{WDPath} for $\operatorname{BES}_h(3)$ (for any $h>0$), which states that it can be decomposed into a pre-infimum process behaving as a Brownian motion, and a post-infimum process behaving as a $\operatorname{BES}(3)$ process started at the infimum. 

\begin{lemma}[\protect{\cite{WDPath}}]  \label{wildeco}
For $h>0$, let $J_3 \sim \operatorname{Unif}(0,h)$, and let $B$ be an independent Brownian motion.
Let $\tau_3$ be the first time $B$ hits $J_3 - h$.
Define $W(t)=h+B(t)$ for $t\in[0,\tau_3]$, and let $(W(\tau_3+t)-W(\tau_3))_{t\ge 0}$ be a $\operatorname{BES}(3)$ process, independent of $B$ and $J_3$.
Then $W$ has the same law as $\operatorname{BES}_h(3)$.
\end{lemma}

The first alternative description gives a similar path decomposition for $Z$.

\begin{prop}\label{inf_path_decom}
The process $Z$ can be constructed as follows.
Let $J$ be a random variable with density $\left(\frac{6 x^2}{z^3} - \frac{4 x^3}{z^4}\right)\mathds{1}_{0<x<z}$, and let $B$ be an independent Brownian motion.
Let $\tau_J$ be the first time $B$ hits $J-z$. Then $Z(t)=z+B(t)$ for $t\in [0, \tau_J]$, and $(Z(\tau_J+t)-Z(\tau_J))_{t\ge 0}$ is a $\operatorname{BES}(3)$ process, independent of $B$ and $J$.
\end{prop}

\begin{proof}
From Proposition~\ref{limZdef} and the Williams path decomposition (Lemma~\ref{wildeco}), the stated decomposition holds if we take $J=\min_{t\ge 0}Z(t)$.
It therefore suffices to verify that the probability density of $\min_{t\ge 0}Z(t)$ is indeed $\left(\frac{6 x^2}{z^3} - \frac{4 x^3}{z^4}\right)\mathds{1}_{0<x<z}$.

For any $0 < p < k/z < 1$, the joint density of $P$ and $\tilde K$ at $(p, k)$ equals
\[
\frac{z}{k} \cdot \frac{6k(z-k)}{z^3}.
\]
With probability $P$, we have $\min_{t\ge 0}Z(t)=\tilde K$. With probability $1-P$, $\min_{t\ge 0}Z(t)$ equals $\tilde K$ plus the minimum of $\operatorname{BES}_{z-\tilde K}(3)$, which (conditional on $P$ and $\tilde K$) is $\operatorname{Unif}(\tilde K, z)$ by Lemma~\ref{infimum}.
Therefore, the density of $\min_{t\ge 0}Z(t)$ at $x$ equals
\[
\iint_{0<p<k/z<1} \frac{1}{z-k}\,(1-p)\,
\frac{z}{k}\,\frac{6 k(z-k)}{z^3}\,dp\,dk
+ \int_0^{x/z} p\,\frac{z}{x}\,\frac{6 x(z-x)}{z^3}\,dp
= \frac{6 x^2}{z^3} - \frac{4 x^3}{z^4}.
\]
Thus, the conclusion follows.
\end{proof}

As a byproduct, this gives another disproof of Conjecture~\ref{coj:main}.

\begin{proof}[Another proof of \Cref{thm:nb5}]
If $X$ were Markovian and (by Proposition~\ref{OP_Prop_infi} and Remark~\ref{OP_Rem_infi}) had the same law as $\operatorname{BES}(5)$, the law of $Z$ would be $\operatorname{BES}_z(5)$.
From Lemma~\ref{infimum}, the density of the infimum of $\operatorname{BES}_z(5)$ is $\frac{3 x^2}{z^3}$, which differs from that of $J$ in Proposition~\ref{inf_path_decom}.
Therefore, we obtain a contradiction.
\end{proof}

From Proposition~\ref{inf_path_decom} and by applying Lemma~\ref{wildeco} again, we obtain another description of $Z$ as a mixture of $\operatorname{BES}_h(3)$ processes for $h \in [0, z]$.

\begin{prop}  \label{BESmix}
The law of $Z$ is $z - H + R_{3,H}$, where $H$ is a random variable taking values in the interval $[0, z]$, with density at $h \in [0, z]$ given by $12h^2(z-h)z^{-4}$, and conditional on $H$, $R_{3,H}$ is a $\operatorname{BES}_H(3)$ process.
\end{prop}

\begin{proof}
For the process $z - H + R_{3,H}$, by Lemma~\ref{infimum}, the density of its minimum at $x$ is
\[
\int_{z-x}^z h^{-1}\,12h^2(z-h)z^{-4}\,dh = \frac{6 x^2}{z^3} - \frac{4 x^3}{z^4}.
\]
Hence, by Proposition~\ref{inf_path_decom} and Lemma~\ref{wildeco}, the conclusion follows.
\end{proof}

\subsection{Future infimum reflection}  \label{s:fir}
Following Pitman's $2M-B$ theorem, it has been shown (first in \cite{Jeulin_enlargement}) that for $2M-B$ (i.e., $\operatorname{BES}(3)$), if it is reflected against its future infimum, one recovers the Brownian motion.
Such a property can also be viewed as the semimartingale decomposition of $\operatorname{BES}(3)$ in an enlarged filtration obtained by including the future infimum.

We next illustrate that a similar future-infimum reflection also holds for $Z$.

\begin{prop}  \label{futinf}
Let $J(t) = \inf_{u \ge t} Z(u)$ be the future infimum process of $Z$, and let $\mathcal{G}_t := \mathcal{F}^{Z}_{t} \lor \sigma(J(t))$ be the sigma-algebra generated by $J(t)$ and $\mathcal{F}^Z_t$. 
Then the process $Z(t) - 2(J(t)-J(0)) - z$ is a $\left(\mathcal{G}_t\right)_{t \ge 0}$-Brownian motion.
\end{prop}

\begin{proof}
Let $\tau$ denote the first time $Z$ hits $J(0)$, and define
\begin{equation}   \label{eq:defr}
R(t) = Z(\tau+t) - Z(\tau), \quad t \ge 0.
\end{equation}
Note that $\tau$ is a stopping time with respect to the filtration $(\mathcal{G}_{t})_{t\ge 0}$. By Proposition~\ref{inf_path_decom}, it is clear that $2J - Z$ on $[0, \tau]$ is a Brownian motion starting at $2J(0) - z$, adapted to $(\mathcal{G}_{t})_{t\ge 0}$, and that $R$ is a $\operatorname{BES}(3)$ process.

Now consider $2J - Z$ on $[\tau, \infty)$.
Since $Z(\tau) = J(\tau)$, we have for any $t \ge 0$,
\begin{equation}   \label{eq:2JZ}
(2J(\tau+t) - Z(\tau + t)) - (2J(\tau) - Z(\tau))
= 2\big(J(\tau+t) - Z(\tau)\big) - \big(Z(\tau + t) - Z(\tau)\big).
\end{equation}
Let $\tilde{J}(t) = \inf_{u \ge t} R(u)$. Then, by \eqref{eq:defr} and the definition of $J$, we have
\[
\tilde{J}(t) = \inf_{u \ge t} Z(\tau+u) - Z(\tau) = J(\tau+t) - Z(\tau).
\]
Therefore, \eqref{eq:2JZ} equals $2\tilde{J}(t) - R(t)$.  
By the future-infimum reflection of $\operatorname{BES}(3)$ (see \cite{Jeulin_enlargement} or \cite[Chapter VI, Corollary 3.7]{revuz_yor}), $2 \tilde{J} - R$ is a Brownian motion adapted to $(\mathcal{G}_{\tau+t})_{t\ge 0}$. Hence, so is $\big((2J(\tau+t)-Z(\tau + t)) - (2J(\tau)-Z(\tau))\big)_{t\ge 0}$.
Therefore, the conclusion follows.
\end{proof}

It is natural to ask whether such a future-infimum reflection holds for $X$ as well. 
However, we expect this not to hold, as this reflection seems to rely on the slope $K'$ being zero.

\subsection{Multi-points distribution}

An implication of Proposition~\ref{OP_Prop_density} is that the distribution of $X(1)$, which is $\operatorname{Chi}(k)$ as pointed out in \cite{OP}, matches the one-point distribution of $\operatorname{BES}(5)$, whereas the multi-point distribution of $X$ appears to be less tractable.  
As an application of Proposition~\ref{BESmix}, we derive the multi-point probability density function of $Z$, which can be viewed as a degeneration of the multi-point probability density function of the original process $X$, thereby shedding light on the latter.

Denote by $p(x,y;t)$ the transition probability density of Brownian motion, namely,
\[
p(x,y;t) = \frac{\phi((x-y)/\sqrt{t})}{\sqrt{t}} = \frac{1}{\sqrt{2\pi t}} \exp\left(-\frac{(x-y)^2}{2t}\right).
\]
For $x, y > 0$, define
\[
p_*(x,y;t) = p(x,y;t) - p(x,-y;t).
\]
Note that for $\operatorname{BES}(3)$, its transition probability density is given by $\frac{y}{x}p_*(x,y;t)$.

\begin{lemma}\label{mulptZ}
Take any $m \in \N$, $0 < t_1 < \cdots < t_m$, and real numbers $x_1, \ldots, x_m > 0$.
The probability density of $(Z(t_1), \ldots, Z(t_m))$ at $(x_1, \ldots, x_m)$ is
\begin{equation}  \label{eq:Ztfin}
\int_{(z-x_*)\vee 0}^z 12(x_m + h - z)h(z-h)z^{-4} 
\prod_{i=1}^m p_*(x_{i-1}+h-z,\, x_i+h-z;\, t_i-t_{i-1})\, dh,
\end{equation}
where (for convenience of notation) we set $t_0 = 0$, $x_0 = z$, and $x_* = \min\{x_1, \ldots, x_m\}$.
\end{lemma}

\begin{proof}
For any $h > 0$, the joint probability density of the process $\operatorname{BES}_h(3)$ at times $0 < t_1 < \cdots < t_m$ and values $x_1, \ldots, x_m > 0$ is
\[
\prod_{i=1}^m \frac{x_i}{x_{i-1}}\, p_*(x_{i-1}, x_i; t_i - t_{i-1})
= \frac{x_m}{x_0} \prod_{i=1}^m p_*(x_{i-1}, x_i; t_i - t_{i-1}),
\]
where $t_0 = 0$ and $x_0 = h$.
Then, by Proposition~\ref{BESmix}, the conclusion follows.
\end{proof}

By evaluating the integral~\eqref{eq:Ztfin}, we obtain a more explicit expression.

\begin{prop}\label{Ztfinev}
The integral~\eqref{eq:Ztfin} equals
\[
\sum_{\Lambda \subset \{1, \ldots, m\}} (-1)^{|\Lambda|} \Xi_\Lambda 
\prod_{i \in \{1, \ldots, m\} \setminus \Lambda}
\frac{\phi((x_{i-1} - x_i)/\sqrt{t_i - t_{i-1}})}{\sqrt{t_i - t_{i-1}}},
\]
where 
\[
\Xi_\emptyset = z^{-4}\bigl(3(x_* \wedge z)^4 - 4(x_* \wedge z)^3(z + x_m) + 6(x_* \wedge z)^2 x_m z\bigr),
\]
and for any $\Lambda \neq \emptyset$,
\begin{align*}
\Xi_\Lambda = &\, 12 a_\Lambda^{-4} c_\Lambda z^{-4} \Big[
\bigl(b_\Lambda^2 + 2 + a_\Lambda b_\Lambda (z + x_m) + a_\Lambda^2 x_m z\bigr)\phi(b_\Lambda)\\
& - \bigl(b_\Lambda^2 + 2 + a_\Lambda b_\Lambda (z + x_m) + a_\Lambda^2 x_m z 
+ a_\Lambda^2(x_* \wedge z)^2 
- a_\Lambda(x_* \wedge z)(b_\Lambda + a_\Lambda z + a_\Lambda x_m)\bigr) \\
&\quad \times \phi(a_\Lambda(x_* \wedge z) + b_\Lambda) \\
& - \bigl(3b_\Lambda + a_\Lambda z + a_\Lambda x_m + b_\Lambda(b_\Lambda + a_\Lambda z)(b_\Lambda + a_\Lambda x_m)\bigr)\\
&\quad \times \frac{1}{2}\bigl(\erf((a_\Lambda(x_* \wedge z) + b_\Lambda)/\sqrt{2}) - \erf(b_\Lambda/\sqrt{2})\bigr)
\Big],
\end{align*}
with
\[
a_\Lambda = 2\sqrt{\sum_{i \in \Lambda} \frac{1}{t_i - t_{i-1}}}, \qquad
b_\Lambda = -\frac{2\sum_{i \in \Lambda} \frac{x_{i-1} + x_i}{t_i - t_{i-1}}}{a_\Lambda}
= -\frac{\sum_{i \in \Lambda} \frac{x_{i-1} + x_i}{t_i - t_{i-1}}}{\sqrt{\sum_{i \in \Lambda} \frac{1}{t_i - t_{i-1}}}},
\]
and
\[
c_\Lambda =
\frac{\exp\left(-\frac{1}{2} \sum_{i \in \Lambda} \frac{(x_{i-1} + x_i)^2}{t_i - t_{i-1}} + \frac{b_\Lambda^2}{2}\right)}
{(\sqrt{2\pi})^{|\Lambda| - 1}\prod_{i \in \Lambda} \sqrt{t_i - t_{i-1}}}
= \frac{\exp\left(-\frac{1}{2\sum_{i \in \Lambda} \frac{1}{t_i - t_{i-1}}}
\sum_{i,j \in \Lambda} \frac{(x_{i-1} + x_i - x_{j-1} - x_j)^2}{(t_i - t_{i-1})(t_j - t_{j-1})}\right)}
{(\sqrt{2\pi})^{|\Lambda| - 1}\prod_{i \in \Lambda} \sqrt{t_i - t_{i-1}}}.
\]
\end{prop}

To prove this, it is convenient to record the following identities for Gaussian integrals.

\begin{lemma}\label{antiphi}
Antiderivatives of $\phi(x)$, $x\phi(x)$, $x^2\phi(x)$, and $x^3\phi(x)$ are, up to an additive constant, 
$\frac{1}{2}\erf(x/\sqrt{2})$, $-\phi(x)$, $-x\phi(x) + \frac{1}{2}\erf(x/\sqrt{2})$, and $-(x^2 + 2)\phi(x)$, respectively.
\end{lemma}

\begin{proof}[Proof of Proposition \ref{Ztfinev}]
By letting $g=z-h$, \eqref{eq:Ztfin} can be written as
\begin{multline*}
\sum_{\Lambda\subset \{1,\ldots,m\}} (-1)^{|\Lambda|} \prod_{i\in\{1,\ldots,m\}\setminus\Lambda} \frac{\phi((x_{i-1}-x_i)/\sqrt{t_i-t_{i-1}})}{\sqrt{t_i-t_{i-1}}}
\\
\times \int_0^{x_*\wedge z} 12(x_m-g)g(z-g)z^{-4}  \prod_{i\in\Lambda} \frac{\phi((x_{i-1}+x_i-2g)/\sqrt{t_i-t_{i-1}})}{\sqrt{t_i-t_{i-1}}}   dg.
\end{multline*}
When $\Lambda=\emptyset$, the second line equals
\[
z^{-4}\bigl(3(x_*\wedge z)^4 - 4(x_*\wedge z)^3(z+x_m)+6(x_*\wedge z)^2x_mz\bigr).
\]
When $\Lambda\neq\emptyset$, the second line can be written as
\[
\int_0^{x_*\wedge z} 12(x_m-g)g(z-g)z^{-4} \, c_\Lambda\,\phi(a_\Lambda g + b_\Lambda)\, dg.
\]
By a change of variables, this becomes
\[
12a_\Lambda^{-4}c_\Lambda z^{-4} \int_{b_\Lambda}^{a_\Lambda(x_*\wedge z) + b_\Lambda}  (g-a_\Lambda x_m-b_\Lambda)(g-b_\Lambda)(g-a_\Lambda z-b_\Lambda) \, \phi(g) \, d g.
\]
By Lemma \ref{antiphi}, the antiderivative of the integrand is
\begin{multline*}
\bigl(-(g^2+2)+(3b_\Lambda + a_\Lambda z + a_\Lambda x_m)g - (3b_\Lambda^2 + 2a_\Lambda b_\Lambda z + 2a_\Lambda b_\Lambda x_m + a_\Lambda^2 x_m z)\bigr)\phi(g)
\\ - \bigl((3b_\Lambda + a_\Lambda z + a_\Lambda x_m)+ b_\Lambda(b_\Lambda + a_\Lambda z)(b_\Lambda + a_\Lambda x_m)\bigr) \cdot \tfrac{1}{2}\erf\!\bigl(g/\sqrt{2}\bigr).
\end{multline*}
Therefore, the integral above equals
\begin{align*}
&(b_\Lambda^2 + 2 + a_\Lambda b_\Lambda (z + x_m) + a_\Lambda^2 x_m z)\phi(b_\Lambda)\\
    -&(b_\Lambda^2 + 2 + a_\Lambda b_\Lambda (z+x_m) + a_\Lambda^2 x_m z + a_\Lambda^2(x_*\wedge z)^2 - a_\Lambda(x_*\wedge z)(b_\Lambda+a_\Lambda z+a_\Lambda x_m)   ) \\ &\times \phi(a_\Lambda(x_*\wedge z)+b_\Lambda)
\\ - &(3b_\Lambda + a_\Lambda z + a_\Lambda x_m + b_\Lambda(b_\Lambda + a_\Lambda z)(b_\Lambda + a_\Lambda x_m))\\ & \times \tfrac{1}{2}\bigl(\erf((a_\Lambda(x_*\wedge z)+b_\Lambda)/\sqrt{2}) - \erf(b_\Lambda/\sqrt{2})\bigr).
\end{align*}
Finally, by summing over $\Lambda$, the conclusion follows.
\end{proof}

Taking $m = 1$, we obtain the one-point distribution density of $Z$.
Note that this corresponds to the $s \to \infty$ limit of the two-point distribution density of $(X(s), X(s + t))$ at $(z, x)$.

\begin{cor}
For any $t, x > 0$, the probability density of $Z(t)$ at $x$ equals
\begin{align*}
&z^{-4} (x \wedge z)^3 (2x + 2z - 3(x \wedge z)) 
\frac{\phi((x - z)/\sqrt{t})}{\sqrt{t}}\\
&\quad - \frac{3}{4} t^{3/2} z^{-4}
\left( 2 - \frac{(x-z)^2}{t} \right) 
\phi((x + z)/\sqrt{t})\\
&\quad + \frac{3}{4} t^{3/2} z^{-4}
\left( 2 - \frac{|x-z|(x+z)}{t} \right) 
\phi((x - z)/\sqrt{t})\\
&\quad - \frac{3}{8} z^{-4} (t - (z - x)^2)(z + x)
\bigl(\erf(-|x - z|/\sqrt{2t}) - \erf(-(x + z)/\sqrt{2t})\bigr).
\end{align*}
\end{cor}

This clearly implies that $Z$ is not distributed as $\operatorname{BES}_z(5)$, providing yet another proof of \Cref{thm:nb5}.

\subsection{Infinitesimal generator}
In Proposition~\ref{OP_Prop_infi}, it is shown that the infinitesimal generator of $X$ coincides with that of $\operatorname{BES}(5)$.  
This property can also be understood from the perspective of a mixture of $\operatorname{BES}(3)$ processes, although from this viewpoint, the matching of infinitesimal generators appears more coincidental than structural.

More precisely, consider a general mixture of $\operatorname{BES}(3)$ processes defined as follows.  
Take any probability measure on $[0, z]$ that is absolutely continuous with respect to the Lebesgue measure, and let $\gamma:[0, z]\to\R_{\ge 0}$ denote its density function, assumed to be bounded and continuous at $z$.  
Let $\tilde Z$ be a random process constructed the same way as $Z$ in Proposition~\ref{inf_path_decom}, except that $\gamma$ serves as the density function of $J$.  
Namely, let $J$ be a random variable with probability density $\gamma$, let $B$ be an independent Brownian motion, and let $\tau_J$ denote the first time $B$ hits $J-z$.  
Define $\tilde Z(t) = z + B(t)$ for $t \in [0, \tau_J]$, and let $(\tilde Z(\tau_J + t) - \tilde Z(\tau_J))_{t \ge 0}$ be a $\operatorname{BES}(3)$ process, independent of both $B$ and $J$.

\begin{prop}\label{p:infgen}
Let $\varphi$ be a twice continuously differentiable function with compact support on $(0, \infty)$. Then
\[
\lim_{t \to 0_+} \frac{1}{t}\,\mathbb{E}[\varphi(\tilde Z(t)) - \varphi(z)]
= \frac{1}{2}\varphi''(z) + \gamma(z)\varphi'(z).
\]
\end{prop}

\begin{proof}
Take a standard Brownian motion $B$.  
From the above construction of $\tilde Z$, we can couple $\tilde Z$ and $B$ so that on $[0, \tau_J]$, $\tilde Z = B + z$, and on $[\tau_J, \infty)$, $\tilde Z \ge B + z$ almost surely.  
Note that
\[
\lim_{t \to 0_+} \frac{1}{t}\mathbb{E}[\varphi(B(t) + z) - \varphi(z)]
= \frac{1}{2}\varphi''(z),
\]
so it remains to show that
\begin{equation}\label{eq:ZBphi}
\lim_{t \to 0_+} \frac{1}{t}\mathbb{E}[\varphi(\tilde Z(t)) - \varphi(B(t) + z)]
= \gamma(z)\varphi'(z).
\end{equation}

For $\tau_J \ge t$, we have
\[
\mathbb{E}[\varphi(\tilde Z(t)) - \varphi(B(t) + z) \mid \tau_J] = 0.
\]
Conditional on $J$ and $\tau_J$ with $\tau_J < t$, the random variables $(t - \tau_J)^{-1/2}(\tilde Z(t) - J)$ and $(t - \tau_J)^{-1/2}(B(t) + z - J)$ are $\operatorname{Chi}(3)$ and standard normal, respectively.  
Since $\varphi$ is twice continuously differentiable with compact support on $(0, \infty)$, for $\tau_J < t$ we have
\[
\left|\mathbb{E}[\varphi(\tilde Z(t)) - \varphi(B(t) + z) \mid J, \tau_J]
- \frac{4}{\sqrt{2\pi}}(t - \tau_J)^{1/2}\varphi'(J)\right|
< C_1 (t - \tau_J),
\]
for some constant $C_1 > 0$ depending only on $\varphi$.  
Note that
\[
\lim_{t \to 0_+} \frac{1}{t}\mathbb{E}[(t - \tau_J) \vee 0]
\le \lim_{t \to 0_+} \mathbb{P}[\tau_J \le t] = 0.
\]
Therefore, the left-hand side of \eqref{eq:ZBphi} equals
\begin{equation}\label{eq:ZBphi2}
\frac{4}{\sqrt{2\pi}}
\lim_{t \to 0_+} \frac{1}{t}
\mathbb{E}\!\left[((t - \tau_J) \vee 0)^{1/2} \varphi'(J)\right].
\end{equation}

We also have that, since $\varphi$ is twice continuously differentiable with compact support on $(0, \infty)$, there exists a constant $C_2 > 0$ depending only on $\varphi$ such that
\begin{equation}\label{eq:ZBphierr}
\lim_{t \to 0_+} \frac{1}{t}\mathbb{E}\!\left[((t - \tau_J) \vee 0)^{1/2} |\varphi'(J) - \varphi'(z)|\right]
< C_2 \lim_{t \to 0_+} \frac{1}{t}\mathbb{E}\!\left[((t - \tau_J) \vee 0)^{1/2} |z - J|\right].
\end{equation}
Note that
\[
\mathbb{E}\!\left[((t - \tau_J) \vee 0)^{1/2} \mid J\right]
= \int_0^t f_{0, z - J}(w)(t - w)^{1/2}\,dw
< C_3 t^{1/2} \exp\!\left(-\frac{(z - J)^2}{2t}\right),
\]
for some universal constant $C_3 > 0$.  
Therefore, since the probability density of $J$ is bounded, the limit on the right-hand side of \eqref{eq:ZBphierr} equals zero.

Now \eqref{eq:ZBphi2}, or equivalently the left-hand side of \eqref{eq:ZBphi}, equals
\[
\varphi'(z) \cdot \frac{4}{\sqrt{2\pi}}
\lim_{t \to 0_+} \frac{1}{t} \mathbb{E}\!\left[((t - \tau_J) \vee 0)^{1/2}\right]
= \varphi'(z) \cdot \frac{4}{\sqrt{2\pi}}
\lim_{t \to 0_+} \frac{1}{t}
\int_0^t \int_0^z \gamma(h) f_{0, z - h}(w) (t - w)^{1/2}\, dh\, dw.
\]
By a change of variables and using the expression for $f_{0, z - h}$, we have
\begin{align*}
&\frac{1}{t}\int_0^t \int_0^z (\gamma(h) - \gamma(z)) f_{0, z - h}(w) (t - w)^{1/2}\, dh\, dw \\
&\quad = \int_0^1 \int_0^{z t^{-1/2}}
(\gamma(z - g t^{1/2}) - \gamma(z))
\cdot \frac{g}{\sqrt{2\pi w^3}}
\exp\!\left(-\frac{g^2}{2w}\right)
(1 - w)^{1/2}\, dg\, dw.
\end{align*}
By the continuity of $\gamma$ at $z$ and the boundedness of $\gamma$, the above integral tends to zero as $t \to 0_+$.  
Hence the left-hand side of \eqref{eq:ZBphi} equals
\[
\varphi'(z)\gamma(z) \cdot \frac{4}{\sqrt{2\pi}}
\lim_{t \to 0_+} \frac{1}{t}
\int_0^t \int_0^z f_{0, z - h}(w) (t - w)^{1/2}\, dh\, dw.
\]
Then we write
\[
\frac{1}{t}\int_0^t \int_0^z  f_{0,z-h}(w) (t-w)^{1/2} \,dh\, dw
=
\frac{1}{t}\int_0^t \int_0^z
\frac{z - h}{\sqrt{2\pi w^3}}
\exp\!\left(-\frac{(z - h)^2}{2w}\right)
(t - w)^{1/2}\, dh\, dw.
\]
After a change of variables, this becomes
\[
\int_0^1 \int_0^{z t^{-1/2}}
\frac{g}{\sqrt{2\pi w^3}}
\exp\!\left(-\frac{g^2}{2w}\right)
(1 - w)^{1/2}\, dg\, dw.
\]
Integrating with respect to $g$, and then substituting $w \mapsto w^2$, we obtain
\[
\int_0^1
\frac{2}{\sqrt{2\pi}}
(1 - w^2)^{1/2}
\left(1 - \exp\!\left(-\frac{z^2}{2t w^2}\right)\right)
dw.
\]
By the dominated convergence theorem, as $t \to 0_+$ this integral converges to $\frac{\sqrt{2\pi}}{4}$.  
Therefore, \eqref{eq:ZBphi} holds.
\end{proof}

\section{Markovian versus non-Markovian projections of $\Psi$}  \label{sec:fil}
\begin{table}[h]
    \centering
    \begin{tabular}{|c|c|}
\hline
Processes & Markovian verse Non-Markovian \\
\hline
$(K(t), K'(t), B(t))_{t\ge 0}$ & Non-Markovian \\
 $(K(t), B(t))_{t\ge 0}$ & Non-Markovian \\
$(K'(t), B(t))_{t\ge 0}$ & Non-Markovian \\

$(D(t)-t)_{t\ge 0}$ plus any non-empty subset of $\{K, B, X\}$ & Non-Markovian \\

$(D(t)-t)_{t\ge 0}$ & Non-Markovian \\
$(D(t)-t, K(t)-B(t))_{t\ge 0}$ & Non-Markovian \\

$(K^{\prime}(t), D(t) -t)_{t \geq 0}$ & time-homogeneous Markovian \\
$(K(t), K'(t), D(t) -t)_{t \geq 0}$ & time-homogeneous Markovian \\
$(K(t)-B(t), K'(t),  D(t) -t)_{t \geq 0}$ & time-homogeneous Markovian \\

$(X(t), K^{\prime}(t), D(t) -t)_{t \ge 0}$ & Non-Markovian \\
 $(B(t), K^{\prime}(t), D(t) -t)_{t \ge 0}$ & Non-Markovian \\
\hline
\end{tabular}
\caption{This table summaries the Markovian and non-Markovian projections.}
\end{table}
So far, we have seen that $X = 2K - B$ is non-Markovian, although it is a projection of the time-homogeneous Markov process $\Psi(t) = (K'(t), K(t), K(t) - B(t), D(t) - t)$, $t \ge 0$, from \Cref{OP_THM_markov}. 
It is therefore an intriguing problem to understand how information and the Markov property are lost under this projection; more precisely, whether some other projected processes of $\Psi$ (in particular, those including $X$) are Markovian. 
We investigate this question in this section.

First, as indicated in \cite{OP}, for any projection of $\Psi$ (of dimension $2$ or higher) to be Markovian, it typically needs to include the distance to the next vertex $(D(t) - t)_{t \ge 0}$. 
Such an analysis was carried out in \cite{OP}, right before Theorem~1.2 there.
\begin{lemma}[\protect{\cite{OP}}]
The processes $(K(t), K'(t), B(t))_{t \ge 0}$, $(K(t), B(t))_{t \ge 0}$, and $(K'(t), B(t))_{t \ge 0}$ are not Markovian.
\end{lemma}
In fact, only $(K(t), K'(t), B(t))_{t \ge 0}$ was analyzed explicitly in \cite{OP}, but the same arguments apply to the other two processes stated here. 
Any other two-dimensional projection of $(K, K', B)$ can be investigated similarly.

Now we consider projections of $\Psi$ that include $(D(t) - t)_{t \ge 0}$ but not $K'$.  
\begin{lemma}
For $\Gamma$ being one of $K$, $B$, $X$, $(K, B)$, $(K, X)$, $(B, X)$, or $(K, B, X)$, the process $(D(t) - t, \Gamma(t))_{t \ge 0}$ is not Markovian.
\end{lemma}
\begin{proof}
Let $(\cG_t)_{t \ge 0}$ be the filtration generated by the process $(D(t) - t,\, \Gamma(t))_{t \ge 0}$.

Fix $t > 0$.  
If $\Gamma$ includes $K$, then $K'(t) \in \cG_t$.  
If $\Gamma$ includes $B$ or $X$, then $K'(t) \le K'(G(t)-)$ and $K'(G(t)-) \in \cG_t$ (where $K'(G(t)-)$ denotes the value of $K'$ in a left neighborhood of $G(t)$).  
In either case, for any $\epsilon > 0$, there exists an event $\cE_\epsilon$, measurable with respect to $\cG_t$, such that $\PP[\cE_\epsilon] > 0$ and $\cE_\epsilon$ implies $K'(t) < \epsilon$.

Next, consider $D(D(t)+) - D(t)$ (where $D(D(t)+)$ denotes the value of $D$ in a right neighborhood of $D(t)$).  
By \Cref{OP_THM_markov}, conditional on $K'(t)$, the random variable $D(D(t)+) - D(t)$ is independent of $\cG_t$.  
Therefore, by \Cref{OP_THM_markov}, for any $a > 0$, we have almost surely
\begin{align*}
&\lim_{\epsilon \to 0} \PP[ D(D(t)+) - D(t) < a \mid (D(t) - t,\, \Gamma(t)),\, \cE_\epsilon] \\
&= \lim_{\epsilon \to 0} \PP[ D(D(t)+) - D(t) < a \mid K'(t) < \epsilon] = 0.
\end{align*}

On the other hand, it follows from Proposition~\ref{OP_Prop_density} and \Cref{OP_THM_markov} that almost surely,
\[
\PP[ D(D(t)+) - D(t) < a \mid (D(t) - t,\, \Gamma(t)) ] > 0.
\]
Thus, for sufficiently small $\epsilon$, and for $(D(t) - t,\, \Gamma(t))$ in a set of positive probability, the conditional law of $D(D(t)+) - D(t)$ given $(D(t) - t,\, \Gamma(t))$ differs from that given $(D(t) - t,\, \Gamma(t))$ and $\cE_\epsilon$.  
Hence, the process $(D(t) - t,\, \Gamma(t))_{t \ge 0}$ is not Markovian.
\end{proof}

For projections of $\Psi$ that include $(D(t) - t)_{t \ge 0}$ but none of $K'$, $K$, $B$, or $X$, the remaining natural candidates for being Markovian are $(D(t) - t)_{t \ge 0}$ and $(D(t) - t, K(t) - B(t))_{t \ge 0}$.  
Unfortunately, neither of these processes is Markovian.
\begin{lemma}
The processes $(D(t) - t)_{t \ge 0}$ and $(D(t) - t, K(t) - B(t))_{t \ge 0}$ are not Markovian.
\end{lemma}
\begin{proof}
Fix $t > 0$.  
As in the previous proof, let $D(D(t)+)$ denote the value of $D$ in a right neighborhood of $D(t)$.  
We claim that $G(t)$ and $D(D(t)+) - D(t)$ are not independent conditional on $D(t) - t$.  
Since $G(t)$ is measurable with respect to the sigma-algebra generated by $(D(s) - s)_{0 \le s \le t}$, this claim implies that the process $(D(t) - t)_{t \ge 0}$ is not Markovian.

We now prove the claim.  
Take any $w > 0$.  
By Proposition~\ref{OP_Prop_density}, as $\epsilon \to 0$, the conditional distribution of $K'(t)$ given $D(t) - t = w$ and $G(t) < \epsilon$ converges to $0$.  
By \Cref{OP_THM_markov}, conditional on $K'(t)$, the random variable $D(D(t)+) - D(t)$ is independent of both $D(t) - t$ and $G(t)$.  
Therefore, by \Cref{OP_THM_markov}, for any $a > 0$ we have
\begin{align*}
&\lim_{\epsilon \to 0} \PP[ D(D(t)+) - D(t) < a \mid D(t) - t = w,\, G(t) < \epsilon] \\ 
&\qquad = \lim_{\epsilon \to 0} \PP[ D(D(t)+) - D(t) < a \mid K'(t) < \epsilon] = 0.
\end{align*}

On the other hand, it follows from Proposition~\ref{OP_Prop_density} and \Cref{OP_THM_markov} that
\[
\PP[ D(D(t)+) - D(t) < a \mid D(t) - t = w ] > 0.
\]
Hence, for sufficiently small $\epsilon$, the conditional law of $G(t)$ given $D(t) - t = w$ and $G(t) < \epsilon$ differs from that given only $D(t) - t = w$, and the claim follows.

Similarly, one can deduce that $G(t)$ and $D(D(t)+) - D(t)$ are not independent conditional on $D(t) - t$ and $K(t) - B(t)$.  
Therefore, the process $(D(t) - t,\, K(t) - B(t))_{t \ge 0}$ is also non-Markovian.
\end{proof}

It remains to consider projections of $\Psi$ that include $(K'(t), D(t) - t)_{t \ge 0}$.  
It is straightforward to verify that some projections of $\Psi$ (while not including $X$) are still Markovian.
\begin{lemma}
The processes $(K'(t), D(t) - t)_{t \ge 0}$, $(K(t), K'(t), D(t) - t)_{t \ge 0}$, and $(K(t) - B(t), K'(t), D(t) - t)_{t \ge 0}$ are time-homogeneous Markov processes.
\end{lemma}
This follows directly from \Cref{Pitman_seq}, and we omit the details.

On the other hand, the Markov property is violated when one of $X$ or $B$ is further included.
\begin{lemma}\label{lem:non-mark}
The processes $(X(t), K'(t), D(t) - t)_{t \ge 0}$ and $(B(t), K'(t), D(t) - t)_{t \ge 0}$ are not Markovian.
\end{lemma}
\begin{proof}
Let $(\cG_t)_{t \ge 0}$ be the filtration generated by $(X(t), K'(t), D(t) - t)_{t \ge 0}$.  
Then $X(D(t)) \in \cG_t$.  
Indeed, since $G(t) \in \cG_t$, we have $X(G(t)) \in \cG_t$, and
\begin{align*}
X(D(t)) &= K(D(t)) = K(G(t)) + (D(t) - G(t))K'(t) \\
&= X(G(t)) + (D(t) - G(t))K'(t).
\end{align*}
On the other hand,
\[
X(D(t)) = K(D(t)) = K(t) + (D(t) - t)K'(t).
\]
By Proposition~\ref{OP_Prop_density}, $K(t)$ remains random conditional on $(X(t), K'(t), D(t) - t)$; hence $X(D(t))$ is also random under this conditioning.  
These observations imply that the conditional law of $X(D(t))$ given $\cG_t$ differs from that given $(X(t), K'(t), D(t) - t)$.  
Therefore, $(X(t), K'(t), D(t) - t)_{t \ge 0}$ is not Markovian.

Similarly, if we let $(\cH_t)_{t \ge 0}$ denote the filtration generated by $(B(t), K'(t), D(t) - t)_{t \ge 0}$, then $B(D(t)) \in \cH_t$.  
However, by Proposition~\ref{OP_Prop_density}, $B(D(t))$ remains random conditional on $(B(t), K'(t), D(t) - t)$.  
Hence $(B(t), K'(t), D(t) - t)_{t \ge 0}$ is also not Markovian.
\end{proof}

The above proof indicates that if one wishes to construct a Markov process that includes $(X(t), K'(t), D(t) - t)_{t \ge 0}$, the only possible choice is $\Psi$.  
More precisely, for such a process to be Markovian, the law of $X(D(t))$ conditional on its value at time $t$ should coincide with the law of $X(D(t))$ conditional on its entire history up to time $t$.  
Note that the latter sigma-algebra contains $\cG_t$ from the preceding proof; therefore, $X(D(t))$ must be determined by the process at time $t$, which implies that $K$ must be included in the process.  
In this way, one is led back to $\Psi$.  
The same reasoning applies to $(B(t), K'(t), D(t) - t)_{t \ge 0}$.

\section{The filtration of $X$}  \label{sec:fX}
In this section, we consider the filtration $(\cF_t^X)_{t \ge 0}$.  
In Pitman's $2M - B$ theorem, the filtrations of $2M - B$, $M - B$, and $B$ satisfy the following relation:
\begin{equation}\label{eq:inclu}
\mathcal{F}_t^{2M - B} \subsetneq \mathcal{F}_t^{M - B} = \mathcal{F}_t^{B}.    
\end{equation}
As seen in \Cref{s:fir}, if we further consider the future infimum of $2M - B$, that is,
$J(t) = \inf_{u \ge t} \bigl(2M(u) - B(u)\bigr)$, then
\[
\mathcal{F}_t^{2M - B} \vee \sigma(J(t)) = \mathcal{F}_t^{M - B} = \mathcal{F}_t^{B}.
\]
For details, see, e.g.,~\cite[Chapter~VI]{revuz_yor}.

An inclusion of filtrations like \eqref{eq:inclu} is not available in the setting of $X = 2K - B$.  
However, if we include either the current location or the slope of the concave majorant over the infinite time horizon (i.e., some information about the Brownian motion up to time infinity) the following relation holds.
\begin{prop}\label{filtering_prop_1}
For any $t \ge 0$,
\[
\mathcal{F}_t^X \vee \sigma(K(t)) \vee \sigma(K'(t))
\subset \mathcal{F}_t^B \vee \mathcal{F}_t^K
= \mathcal{F}_t^{K - B} \vee \mathcal{F}_t^{K'}
= \mathcal{F}_t^{B} \vee \sigma(K(t))
= \mathcal{F}_t^{B} \vee \sigma(K'(t)).
\]
\end{prop}

\begin{proof}
The first inclusion is immediate.

We now show that all the sigma-algebras in the stated equality are equivalent to $\mathcal{F}_t^B \vee \mathcal{F}_t^K$, which represents the information of the Brownian path up to time~$t$ together with the concave majorant over the infinite time horizon, restricted to~$[0,t]$.

For the first equality, note that $\mathcal{F}_t^{K'} = \mathcal{F}_t^K$, so
\[
\mathcal{F}_t^B \vee \mathcal{F}_t^K
= \mathcal{F}_t^{K - B} \vee \mathcal{F}_t^K
= \mathcal{F}_t^{K - B} \vee \mathcal{F}_t^{K'}.
\]

For $\mathcal{F}_t^{B} \vee \sigma(K'(t))$, the inclusion 
$\mathcal{F}_t^{B} \vee \sigma(K'(t)) \subset \mathcal{F}_t^B \vee \mathcal{F}_t^K$
is obvious.  
Conversely, one can reconstruct $K$ on~$[0,t]$ as follows:
construct the concave majorant of $B$ on~$[0,t]$, remove all line segments whose slope is smaller than~$K'(t)$, and then attach a single line segment with slope~$K'(t)$.  
This procedure recovers~$K$ on~$[0,t]$, and therefore
$\mathcal{F}_t^B \vee \mathcal{F}_t^K \subset \mathcal{F}_t^{B} \vee \sigma(K'(t))$.

A similar argument applies to $\mathcal{F}_t^{B} \vee \sigma(K(t))$.  
The inclusion $\mathcal{F}_t^{B} \vee \sigma(K(t)) \subset \mathcal{F}_t^B \vee \mathcal{F}_t^K$ is clear.  
Conversely, one can reconstruct~$K$ on~$[0,t]$ by first forming the concave majorant of $B$ on~$[0,t]$, then deleting any line segment whose intercept at time~$t$ equals~$K(t)$, and finally adding a line segment connecting the right end of the last remaining piece to the point~$(t, K(t))$.  
This yields $K$ on~$[0,t]$, and thus
$\mathcal{F}_t^{B} \vee \sigma(K(t)) = \mathcal{F}_t^B \vee \mathcal{F}_t^K$.
\end{proof}

In the case of $2M - B$, the future infimum reflection (from~\cite{Jeulin_enlargement}; see \Cref{s:fir}) can be interpreted as a semimartingale decomposition of $2M - B$ (i.e., $\operatorname{BES}(3)$) under an enlarged filtration obtained by including the future infimum.  
In analogy, a similar enlargement of filtration arises for $X = 2K - B$, where the additional information consists of $(K(t) - B(t))$ and $(D(t) - t)$.  
To the best of our knowledge, this constitutes a new type of progressive enlargement of filtration, since existing examples in the literature involve only random times or future infimum.

\begin{cor}
Define 
\[
\mathcal{G}_t := \mathcal{F}_t^{B} \vee \sigma(K^{\prime}(t)) \vee \sigma(D(t) - t).
\]
Then $(\mathcal{G}_t)_{t \ge 0}$ is a filtration.  
Moreover, both $B$ and $X$ are semimartingales with respect to this filtration, and they admit the following decompositions:
\[
B(t) = \tilde{B}(t) + \int_0^t \left( K^{\prime}(s) + \frac{1}{K(s) - B(s)} - \frac{K(s) - B(s)}{D(s) - s} \right) ds,
\]
\[
X(t) = -\tilde{B}(t) + \int_0^t \left( K^{\prime}(s) - \frac{1}{K(s) - B(s)} + \frac{K(s) - B(s)}{D(s) - s} \right) ds,
\]
where $\tilde{B}$ is a Brownian motion under the filtration $(\mathcal{G}_t)_{t \ge 0}$.
\end{cor}

\begin{proof}
By Proposition~\ref{filtering_prop_1}, we have 
\[
\mathcal{G}_t = \mathcal{F}_t^{B} \vee \mathcal{F}_t^{K} \vee \sigma(D(t) - t),
\]
which is clearly a filtration.  
From \Cref{OP_THM_markov}, conditional on $\mathcal{G}_t$, the process $X - K = K - B$ is a 3-dimensional Bessel bridge from $(t, K(t) - B(t))$ to $(D(t), 0)$.  
The stated decompositions for $B$ and $X = 2K - B$ then follow from the stochastic differential equation satisfied by this bridge (Lemma \ref{lem:sdebb}).
\end{proof}

\section{Simulations}  \label{sec:simu}
\begin{figure}[!ht]
    \centering
\begin{subfigure}{0.6\hsize}
    \resizebox{0.99\hsize}{!}{
\includegraphics[]{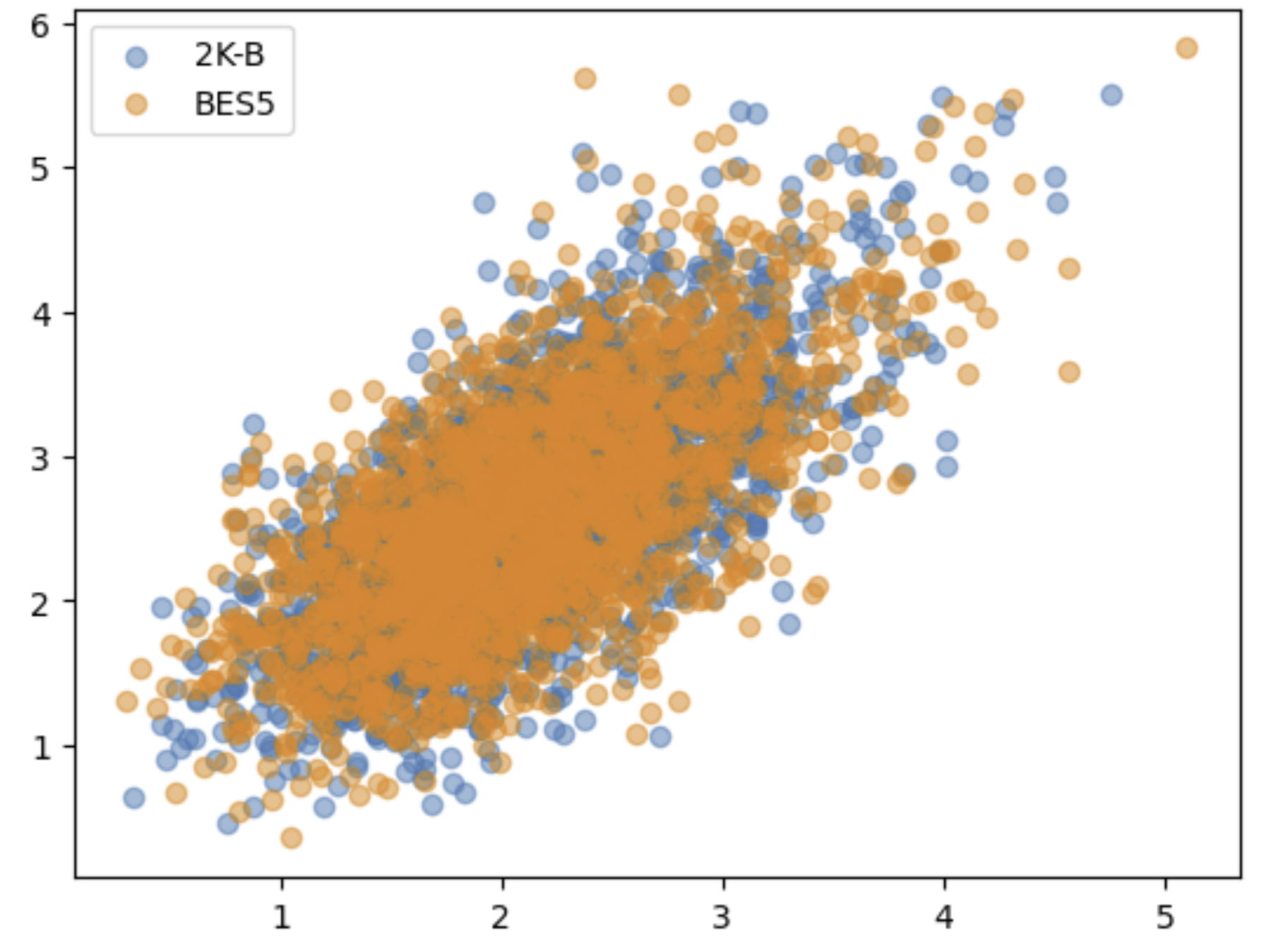}
    
}
\end{subfigure}
    \caption{Scatter plot for 3000 samples of $(X(1),X(2))$ and $(Y(1),Y(2))$}
    \label{scatter_plot}
\begin{subfigure}{0.6\hsize}
    \resizebox{0.99\hsize}{!}{
\includegraphics[]{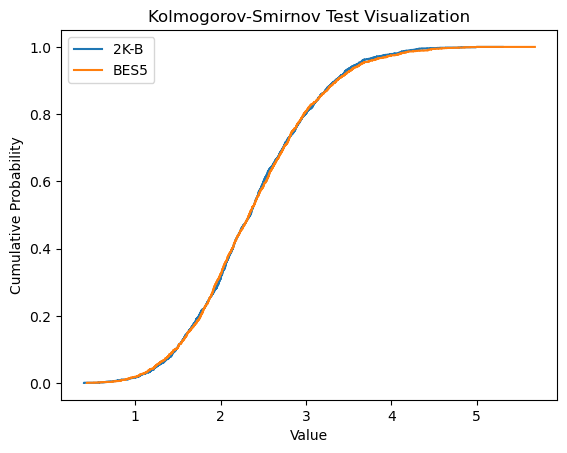}
}
\end{subfigure}
    \caption{KS test for $\lambda_1X(1)+ \lambda_2X(2)$ and $\lambda_1 Y(1)+ \lambda_2 Y(2)$, with $\vec\lambda = (1,1)$, for $3000$ samples.}
    \label{fig:ks2}
\end{figure}

\begin{figure}[!ht]
    \centering

\end{figure}

In this final section, we present numerical simulations comparing $X=2K-B$ and $Y= \operatorname{BES}(5)$.

It is a nontrivial task to simulate the concave majorant on $[0, \infty)$ numerically, since it involves information about the Brownian motion up to infinite time.  
However, this becomes feasible thanks to Theorem~\ref{Pitman_seq} (see also~\cite{carolan2003characterization}).  
More precisely, we can sample $K$ on a finite interval $[t, T]$ using the following recursive procedure.

Let $t_1 := t$, and sample $(K(t_1), K^{\prime}(t_1), K(t_1) - B(t_1), D(t_1) - t_1)$ according to Proposition~\ref{OP_Prop_density}, which determines $K$ on $[t_1, D(t_1)]$.  
We then proceed recursively according to Theorem~\ref{Pitman_seq}.  
For any $n \ge 1$, suppose that we already have $t_n$, $K(t_n)$, $K^{\prime}(t_n+)$, and $D(t_n+)$; we then sequentially take:
\begin{itemize}
    \item $t_{n+1} = D(t_n+)$,
    \item $K(t_{n+1}) = K(t_n) + (t_{n+1} - t_n) K^{\prime}(t_n+)$,
    \item $K^{\prime}(t_{n+1}+) \sim \operatorname{Unif}(0, K^{\prime}(t_n+))$,
    \item $D(t_{n+1}+) - t_{n+1} \sim \operatorname{Gamma}\!\left(\frac{1}{2}, \frac{1}{2} (K^{\prime}(t_{n+1}+))^2\right)$.
\end{itemize}
We repeat this procedure until $t_N \ge T$ (for some $N \ge 1$), thereby obtaining $K$ on $[t, t_N]$.  
Given $K$ on $[t, t_N]$, we then sample $K - B$ on $[t_1, t_2]$ as a 3-dimensional Bessel bridge (by Theorem~\ref{OP_THM_markov}), and on each interval $[t_n, t_{n+1}]$ for $n=2, 3, \ldots, N-1$ as a Brownian excursion (by Theorem~\ref{Pitman_seq}).

We implemented the simulation using the \texttt{stochastic} Python package~\cite{stochastic_py}.

As shown in~\cite{OP}, the distribution of $X = 2K - B$ at each individual time coincides with that of $Y = \operatorname{BES}(5)$.  
A natural next step, therefore, is to compare their two-point distributions, for instance $(X(1), X(2))$ and $(Y(1), Y(2))$.  
However, distinguishing these distributions numerically turns out to be challenging in practice.  
Figure~\ref{scatter_plot} displays 3000 independent samples of $(X(1), X(2))$ and $(Y(1), Y(2))$.

\begin{figure}[!ht]
    \centering
\begin{subfigure}{0.48\hsize}
    \resizebox{0.99\hsize}{!}{
\includegraphics[]{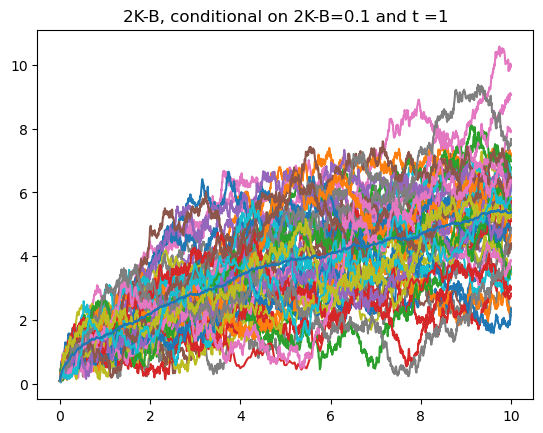}
    
}
\end{subfigure}
\begin{subfigure}{0.48\hsize}
    \resizebox{0.99\hsize}{!}{
\includegraphics[]{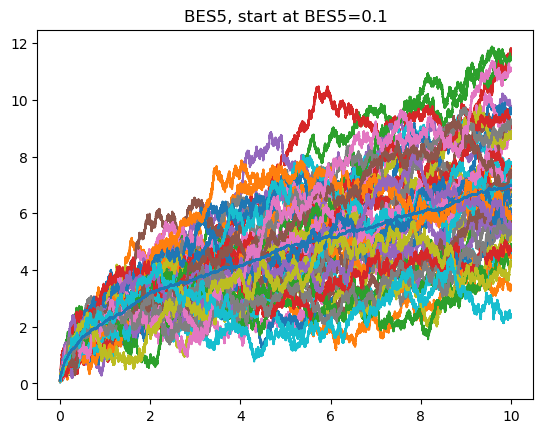}
}
\end{subfigure}
    \caption{50 samples of $X=2K-B$ (left) and $Y= \operatorname{BES}(5)$ (right) conditional on $X(1)=0.1$ and $Y(1)=0.1$, respectively. Here the horizontal axis denotes $t-1$. The dark blue curve is the mean of the simulated trajectories.}
    \label{fig:2k-bconditional}

\begin{subfigure}{0.6\hsize}
    \resizebox{0.99\hsize}{!}{
\includegraphics[]{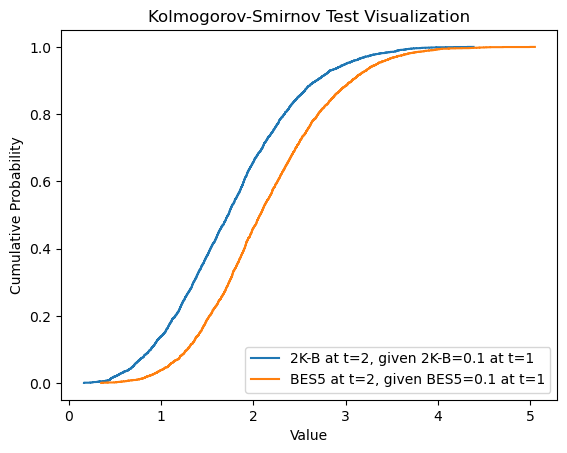}
}
\end{subfigure}
    \caption{KS test for $X(2)$ conditional on $X(1)=0.1$, and $Y(2)$ conditional on $Y(1)=0.1$, with $3000$ samples. The p-value is $2.65 \times10^{-58}$.}
    \label{fig:conditional_KS}
    
\end{figure}
We attempted to distinguish the distributions of $(X(1), X(2))$ and $(Y(1), Y(2))$ using statistical tests.  
Since the Kolmogorov-Smirnov (KS) test applies only in one dimension, we considered one-dimensional projections of these two-dimensional random vectors.  
For a given vector $\vec\lambda = (\lambda_1, \lambda_2)$, we tested the distributions of $\lambda_1 X(1) + \lambda_2 X(2)$ and $\lambda_1 Y(1) + \lambda_2 Y(2)$.  
(Note that any probability distribution on $\mathbb{R}^2$ is determined by all of its one-dimensional projections, via the Fourier transform.)  
We performed KS tests on several such projections with up to 3000 samples, but were unable to reject the null hypothesis that the two random vectors come from the same distribution.  
Figure~\ref{fig:ks2} shows an example of such a test with 3000 samples and $\vec\lambda = (1, 1)$.  
This outcome is not surprising, as the two processes are indeed very similar.

Motivated by our theoretical findings (in particular the first proof of Theorem \ref{thm:nb5}), we adopt a conditional sampling strategy that clearly distinguishes $X$ and $Y$.  
Specifically, we sample $X$ on $[1, \infty)$ conditional on $X(1) = 0.1$, and $Y$ on $[1, \infty)$ conditional on $Y(1) = 0.1$.  
The sampling algorithm described above can be easily adapted to implement this conditional sampling.  
We observe that the process $X$, conditional on $X(1) = 0.1$, exhibits a noticeably weaker upward trend compared with $Y$, conditional on $Y(1) = 0.1$, consistent with Theorem~\ref{thm:nb5} (see Figure~\ref{fig:2k-bconditional}).  
Furthermore, a KS test comparing $X(2)$ conditional on $X(1) = 0.1$ and $Y(2)$ conditional on $Y(1) = 0.1$, each based on $3000$ samples, yields a $p$-value of $2.65 \times 10^{-58}$, strongly rejecting the null hypothesis that the two distributions are identical (see Figure~\ref{fig:conditional_KS}).

\subsection*{Acknowledgement}
The research of LZ was partially supported by the NSF grant DMS-2246664 and by the Miller Institute for Basic Research in Science at the University of California, Berkeley.  
The authors would like to thank Jim Pitman for insightful discussions and for encouraging them to record the discoveries presented in this paper.  
YC also thanks Steven Evans and Alexander Tsigler for valuable discussions, and LZ thanks Dor Elboim and B{\'a}lint Vir{\'a}g, whose talk at the Institute for Advanced Study brought the conjecture in~\cite{OP} to his attention.  
We are further grateful to Alan Hammond for advice on the revision of the first draft, and to the anonymous referees for their careful reading of the manuscript and many valuable suggestions that improved the exposition.

\bibliographystyle{plain}
\bibliography{bibliography}

\end{document}